\documentclass[11pt,leqno]{article}

\usepackage[english]{babel}

\usepackage{comment}
\usepackage{calligra}
\usepackage{mathrsfs}
\usepackage{pb-diagram} 
\usepackage{tikz}
\usepackage{lipsum}
\usepackage{adjustbox}
\usetikzlibrary{cd}
\usetikzlibrary{positioning}
\usepackage{graphics}
\usepackage{mathtools}
\usetikzlibrary{intersections,calc,arrows.meta}
\usepackage{amsmath}
\usepackage{amssymb}
\usepackage{amsfonts}
\usepackage{amsthm}
\usepackage{mathabx}

\newtheorem{theo}{Theorem}[section]%[chapter]

\newtheorem{defi}[theo]{Definition}%[chapter]

\newtheorem{cor}[theo]{Corollary}%[chapter]

\newtheorem{rem}[theo]{Remark}%[chapter]

\newtheorem{prop}[theo]{Proposition}%[chapter]

\newtheorem{lem}[theo]{Lemma}%[chapter]

\newtheorem{ex}[theo]{Example}%[chapter]

\newtheorem{nota}[theo]{Notations}%[chapter]

\newtheorem{conv}[theo]{Convention}
\newtheorem{claim}[theo]{Claim}

\DeclareMathOperator{\Hom}{Hom}
\DeclareMathOperator{\supp}{supp}

\DeclareMathOperator{\intHom}{\mathscr{H}\text{\kern -3pt {\calligra\large om}}\,}
\DeclareMathOperator{\RintHom}{\textit{R}\kern -0pt \mathscr{H}\text{\kern -3pt {\calligra\large om}}\,}
\DeclareMathOperator{\RkintHom}{\textit{R}^{k} \kern -2pt \mathscr{H}\text{\kern -3pt {\calligra\large om}}\,}

\DeclareSymbolFont{largesymbolsA}{U}{txexa}{m}{n}
\DeclareMathSymbol{\varprod}{\mathop}{largesymbolsA}{16}

\newcommand{\mS}{\text{SS}}

\tikzset{
    labl/.style={anchor=south, rotate=90, inner sep=.5mm}
}
\tikzset{
  symbol/.style={
    draw=none,
    every to/.append style={
      edge node={node [sloped, allow upside down, auto=false]{$#1$}}}
  }
}
\newcommand{\function}[5]{%
  \begin{tikzcd}[
    column sep=2em,
    row sep=1ex,
    ampersand replacement=\&
  ]
  #1\colon \&[-3em]
  #2\vphantom{#3} \arrow[r] \&
  #3\vphantom{#2} \\
  \&
  #4\vphantom{#5} \arrow[u,symbol=\in] \arrow[r,mapsto] \&
  #5\vphantom{#4} \arrow[u,symbol=\in]
  \end{tikzcd}%
}
\newcommand{\idfunction}[5]{%
  \begin{tikzcd}[
    column sep=2em,
    row sep=1ex,
    ampersand replacement=\&
  ]
  #1\vphantom{#2} \arrow["#5",r] \&
  #2\vphantom{#1} \arrow[l]\\
  #3\vphantom{#2} \arrow[u,symbol=\in] \arrow[r] \&
  #4\vphantom{#3} \arrow[u,symbol=\in] \arrow[l]
  \end{tikzcd}%
}
\title{The stalk formula for multi-microlocal Hom functors and multi-microlocal Sato’s triangle}
\author{
Ryosuke \textsc{Sakamoto}
}
\date{}

\usepackage[english]{babel}
\begin{document}
\maketitle
\begin{abstract}
The concept of ``multi-microlocalization'' was introduced  to extend the usual microlocal sheaf theory to a more general scope. This paper aims to further extend this theory by exploring advanced topics. One is a stalk formula for multi-microlocalized Hom functors and we compute  some examples in multi-microlocal settings. Secondly we construct the Sato's triangle in the context of multi-microlocal analysis. Ultimately we get the Sato's triangle for the multi-microlocalized Hom functors.
\end{abstract}

\section*{Introduction.}
In this paper, we consider extending the usual microlocal sheaf theory, which holds in the more general scope of the multi-microlocal sheaf theory and subanalytic sheaf theory.  

Sato, Kawai and Kashiwara introduced a new method of analysis, algebraic analysis(\cite{SKK}). This theory concerns the microlocal properties of a function. Kashiwara and Schapira used a homological-geometric method based on the derived category of sheaves to reformalize the theory(\cite{KS85},\cite{KS90}). They introduced the notion of the microsupport $\mS(F)$ for a sheaf $F$ on a derived category of sheaves. This is a formalization of the singularity of a function using the theory of sheaves, and it is the collection of directions in which the cohomology of the sheaf $F$ varies. This work provided a powerful new framework for studying singularities of solutions to differential equations.
The microlocal sheaf theory can also be interpreted as ``a construction that enables the treatment of Fourier transforms on sheaves on manifolds." It is well known that Fourier transforms have applications in differential equation theory. By using the Fourier transformation theory of sheaves, we can consider microlocalization $\mu_M$. Microlocalization is the process of considering a sheaf $F$ on $X$ as a sheaf on the conormal bundle to a submanifold $M$. Based on such foundations, study of differential equations has often been successful.

On the other hand, important analytical objects such as asymptotically developable functions and tempered distributions do not form classical sheaves, therefore conventional methods cannot be applied. New concepts such as Ind-sheaves and sheaves on subanalytic sites have been introduced and studied to handle these objects(\cite{KS01},\cite{Lu1}).  Honda,  Yamazaki, and Prelli (\cite{HPY}) introduced a more extended framework called ``multi-microlocalization." After introducing this multi-microlocalization, they further obtained multi-microlocalization along the universal family of submanifolds of $\Hom$ functors. Algebraic analysis is becoming applicable to differential equations with more complex singularities of a wider class of functions. Various types of microlocalization have been defined and applied to asymptotic analysis and boundary value problems(\cite{T}, \cite{TS}). The theory of multi-microlocalization is constructed to unify these theories.

In Section \ref{sec:pre}, we assume knowledge of the derived category of sheaves and summarize the basic concepts and definitions related to  sheaves  on subanalytic sites and multi-microlocal analysis as explained in \cite{Lu1}, \cite{Lu2} and \cite{HPY}.  Also we would like to give supplementary facts which follow immediately from the past results but could not be found in the literature. These subsections provide a foundation for the rest of the paper, which builds upon this material to explore more advanced topics and ideas. 

In Section \ref{sec:stalkformula}, we derive a stalk formula for multi-microlocalized $\Hom$ functors in general multi-micolocal cases for subanalytic sheaves. A significant challenge in obtaining a stalk formula for subanalytic sheaves is the absence of a microlocal cut-off functor. We circumvent this issue by considering constructible sheaves as the first argument of the $\Hom$-functors. In the following section we prepare some lemmas to obtain more explicit formula for an easy computation in  particular multi-microlcal settings.  Finally we compute some examples in specific multi-microlocal settings. As is well known in the classical theory, a stalk formula is crucial for the microlocal analysis of sheaves (see \cite{KS90}). 
%For example, in \cite{Aoki} a stalk formula is important for analyzing microlocal operators.  
We believe that the stalk formula for multi-microlocalized $\Hom$ functors proved in this paper will be instrumental in the future.
%, the theory multi-microlocal operators in \cite{H} for instance. 

In Section \ref{sec:satotr}, our objective is to establish Sato's triangle within the framework of multi-microlocal analysis. To achieve this, we first prepare Uchida's triangle in subanalytic sheaf theory. Next we need to multi-microlocalize Uchida's triangle in subanalytic sheaf theory. The next step involves multi-microlocalizing Uchida's triangle, which poses a challenge because proper direct image functors do not generally map conic subanalytic sheaves to conic subanalytic sheaves. We address this difficulty by proving several lemmas, ultimately leading to the formulation of Sato's triangle for the $\mu hom$ functor in multi-microlocalizations. We conclude this section with a simple example in a multi-microlocal setting to illustrate the theoretical developments. Experts recognize Uchida's triangle and Sato's triangle as  fundamental tools for the microlocal analysis of sheaves. They are particularly suited for solving microlcal Cauchy problems, boundary value problems and various other issues. Our results extend these concepts within the theory of multi-microlocalizations, opening up possibilities for addressing problems in more general settings.

At the end of the introduction, we would like to show our greatest appreciation to Professor Naofumi Honda for the valuable advice and generous support
in Hokkaido University.

\section{Preliminaries.}\label{sec:pre}
In this section we recall some basic notions. References are made to \cite{Lu1} , \cite{Lu2} and \cite{KS01} for the theory of sheaves on subanalytic sites, and \cite{HPY}, \cite{HP} and \cite{H} for the theory of multi-microlocalizations.
\subsection{Theory of subanalytic sheaves.}
$X$ will be a real analytic manifold and $k$ a field. We assume all manifolds are analytic, finite-dimensional and countable at infinity. Denote by $\mbox{Op}(X_{sa})$(resp. $\mbox{Op}^c(X_{sa})$) the category of open (resp. relatively compact open) subanalytic subsets of $X$. One endows $\mbox{Op}(X_{sa})$ with the following topology: $\mbox{Cov}(U) \subset \mbox{Op}(X_{sa})$ is a covering of $U \in \mbox{Op}(X_{sa})$ if for any compact subset $K$ of $X$ there exists a finite subset $S_0$ of $\mbox{Cov}(U)$ such that $K\cap \bigcup_{V\in S_0} V = K \cap U$. We will call $X_{sa}$ the subanalytic site. Let $\mbox{Mod}(k_{X_{sa}})$ denote the category of sheaves on $X_{sa}$.

Let $\text{Mod}_{\mathbb{R}\text{-}c}(k_X)$ be the abelian category of $\mathbb{R}$-constructible sheaves on $X$, and consider its subcategory $\mbox{Mod}^c_{\mathbb{R}\text{-}c}(k_X)$ consisting of sheaves whose support is compact.
Denote $\rho: X \to X_{sa}$ the natural morphism of sites. $\rho_\ast: \mbox{Mod}(k_{X_{}})\to \mbox{Mod}(k_{X_{sa}})$ is fully faithful and left exact.

\begin{nota}
Since the functor $\rho_\ast$ is fully faithful and exact on $\mbox{Mod}_{\mathbb{R}\text{-}c}(k_X)$ we identify $\mbox{Mod}_{\mathbb{R}\text{-}c}(k_X)$ with the image in $\mbox{Mod}(k_{X_{sa}})$. When there is no risk of confusion we will write $F$ instead of $\rho_\ast F$, for $F \in \mbox{Mod}_{\mathbb{R}\text{-}c}(k_X)$.
\end{nota}
Let $X,Y$ be two real analytic manifolds, and let $f: X\to Y$ be a real analytic map.
We get external operations $f_\ast$ and $f^{-1}$ which are always defined for sheaves on Grothendieck topologies.
Let $Z$ be a locally closed subanalytic set. As in classical sheaf theory we may define $\Gamma_Z$ and $(\cdot)_Z$.
We may define the functor of proper direct image $f_{!!}$ for subanalytic sheaves.
If $f$ is proper on $ \mbox{supp}(F)$ then $f_\ast F \simeq f_{!!}F$, in this case $f_{!!}$ commutes with $\rho_\ast$.

We denote by $D(k_{X_{sa}})$ the derived category of $\mbox{Mod}(k_{X_{sa}})$ and its full subcategory consisting of bounded(resp. bounded below, resp. bounded above) complexes is denoted by $D^b(k_{X_{sa}})$(resp. $D^+(k_{X_{sa}})$, resp. $D^-(k_{X_{sa}})$). As usual we denote by $D^b_{\mathbb{R}\text{-}c}(k_X)$ (resp. $D^b_{\mathbb{R}\text{-}c}(k_{X_{sa}})$) the full subcategory of $D^b(k_X)$ (resp. $D^b(k_{X_{sa}})$) consisting of objects with $\mathbb{R}$-constructible cohomology. We would like to state the result which immediately follows from \cite[Proposition 2.4.8]{Lu1}.
\begin{prop}[]\label{embiso}
\textit{
Let $F \in D^+(k_{Y_{sa}})$, and let $f: X \to Y$ be a closed embedding. Then $Rf_{!!}f^{-1}F\simeq F_X$.
}
\end{prop}
\proof
Since by Proposition Proposition 2.4.8 we have $R\Gamma_X(F) \simeq \RintHom(k_X, F) \simeq Rf_{\ast}f^!F$.
Since $(\cdot)_X$ is a left adjoint to $R\Gamma_X$ i.e. $\Hom(G_X, F) \simeq \Hom(G, R\Gamma_X(F))$, by Yoneda's lemma the result follows.
\endproof
\begin{rem}
Unlike the classical case, this isomorphism does not hold for an open embedding in general.
\end{rem}
The functor $Rf_{!!}$ admits a right adjoint, denoted by $f^!$.
A similar consequence as in \cite[Proposition 5.3.8]{KS01} follows as shown below.
\begin{prop}\label{homiso}
\textit{
Let $G \in D^b(k_{X_{sa}})$, $ F \in D^+(k_{X_{sa}})$. 
There is a natural isomorphism
$
f^!\RintHom(G, F) \simeq \RintHom(f^{-1}G, f^!F).
$
}
\end{prop}
\subsection{Multi-microlocalization.}\label{Sec: multimicro}
Let $\tau_i: E_i \to Z\  (1 \leq i \leq \ell)$ be vector bundles over $Z$, and let $E^\ast_i$ be the dual bundle of $E_i$. We denote by $\wedge_i$ and $\lor_i$ the Fourier-Sato transformations and the inverse Fourier-Sato transformations on $E_i$ respectively. Moreover we denote $\wedge^\ast_i$ and $\lor^\ast_i$ the Fourier-Sato and the inverse Fourier-Sato transformation on $E^\ast_i$ respectively. 
Set $E:= E_1 \underset{Z}{\times}  \cdots \underset{Z}{\times} E_\ell$ and $E^\ast:= E^\ast_1 \underset{Z}{\times} \cdots \underset{Z}{\times} E^\ast_\ell$. 
Let $\tau: E \to Z$ be the canonical projection. Set $P^{\prime}_i := \{(\eta, \xi)\in E_i \underset{Z}{\times} {E_i}^\ast; \langle \eta, \xi \rangle \leq 0\}$. 
\begin{align*}
P^{\prime}:= P^{\prime}_1 \underset{Z}{\times} \cdots \underset{Z}{\times} P^{\prime}_\ell,P^+:= E \underset{Z}{\times} E^\ast\setminus P^{\prime}
\end{align*}
and denote by $p^{\prime}_1: P^{\prime} \to E$ and $p^{\prime}_2: P^{\prime} \to E^\ast$ the canonical projections respectively.
Set $P_i := \{(\eta, \xi)\in E_i \underset{Z}{\times} E^\ast_i; \langle \eta, \xi \rangle \geq 0\}$. 
\[P:= P_1 \underset{Z}{\times} \cdots \underset{Z}{\times} P_\ell,P^-:= E \underset{Z}{\times} E^\ast\setminus P\]
and denote by $p^{\prime\prime}_1: P \to E$ and $p^{\prime\prime}_2: P \to E^\ast$ the canonical projections respectively.
Let $F$ and $G$ be multi-conic objects on $E$ and $E^\ast$, respectively. Then we set for short $\wedge_E$ (resp. $\lor_E^\ast$) the compostion of the Fourier-Sato transformations $\wedge_i$ on $E_i$ (resp. the composition of the inverse Fourier-Sato transformations $\lor_i^\ast$) for each $i \in \{1,\dots, \ell\}$.

Following \cite{HPY}, we recall the notion of multi-microlocalization of sheaves. Let $X$ be a real analytic manifold with $\dim X = n$, and let $\chi= \{M_1, \dots, M_\ell\}$ be a family of closed submanifolds in $X$.
We set, for $N \in \chi$ and $p\in N$,
\[\mbox{NR}_p(N):=\{M_j\in \chi; p\in M_j, N\not \subseteq M_j \mbox{ and }M_j \not \subseteq N\}.\]
Let us consider the following conditions for $\chi$.
\begin{enumerate}
	\item[H1] Each $M_j \in \chi$ is connected and the submanifolds are mutually distinct, i.e., $ M_j \not = M_{j^\prime}$ for $j\not=j^\prime$.
	\item[H2] For any $N \in \chi$ and $p \in N$ with $\mbox{NR}_p(N)\not = \emptyset$, we have
	\begin{align*}
	\left(\bigcap_{M_j \in NR_p(N)}T_p M_j\right) + T_pN = T_pX.
	\end{align*}
\end{enumerate}
If $\chi$ satisfies the condition H2, then for any $p\in X$, there exist a neighborhood $V$ of $p$ in $X$, a system of local coordinates $(x_1, \dots, x_n)$ in $V$ and  a family of subsets $\{I_j\}_{j=1}^\ell$ of the set $\{1,2,\dots, n\}$ for which the following conditions hold.
\begin{enumerate}
    \item[1.] Either $I_k \subset I_j$, $I_j \subset I_k$ or $I_k \cap I_j =\emptyset$ holds $(k,j \in \{1,2,\dots, \ell\})$.
    \item[2.] A submanifold $M_j \in \chi$ with $p \in M_j(j= 1,2, \dots, \ell)$ is defined by $\{x_i =0; i \in I_j\}$ in $V$.
\end{enumerate}
We set, for $N \in \chi$,
\begin{align*}
\iota_\chi(N) := \bigcap_{N\subsetneqq M_j} M_j.
\end{align*}
Here $\iota_\chi(N):= X$ for convention if there exists no $j$ with $N \subsetneqq M_j$. When there is no risk of confusion, we write for short $\iota(N)$.
We also assume the condition
\begin{enumerate}
    \item[H3] $M_j \not = \iota(M_j)$ for any $j \in \{1,2,\dots,\ell\}$.
\end{enumerate}
In local coordinates let $I_1$, $\dots$, $I_\ell$ $\subset \{1,\dots, \ell\}$ such that $M_i = \{x_k=0; k \in I_i\}$. Note that the family $\chi$ satisfies the conditions H1, H2, and H3 if and only if $I_1, \dots, I_\ell$ satisfy the corresponding conditions
\begin{align}\label{1.3}
{}&\mbox{(i) either }I_j \subsetneq I_k, I_k \subsetneq I_j\mbox{ or }I_j \cap I_k = \emptyset\mbox{ holds for any }j \not = k,\\
{}&\mbox{(ii)}\displaystyle\left(\bigcup_{ I_k \subsetneq I_j}I_k\right) \subsetneq I_j\mbox{ for any }j\nonumber.
\end{align}
Hence, for any $j \in \{1,2, \dots, \ell\}$, the set
\begin{align}\label{ihat}
\hat{I}_j:= I_j \setminus \displaystyle\left(\bigcup_{ I_k \subsetneq I_j}I_k\right)
\end{align}
is not empty by the condition H3. By the conditions H1,H2 and H3, we have
for $j \in \{1,\dots,\ell\}$
\[
I_j = \displaystyle\bigsqcup_{I_k \subset I_j} \hat{I}_k.
\]
In particular $ \bigcup_{1\leq j \leq \ell} I_j = \hat{I}_1 \sqcup \dots \sqcup \hat{I}_\ell.$
Set for convenience
$ I_0 = \hat{I}_0 := \{1,\dots, n\} \setminus \bigcup_{j=1}^\ell I_j.$
Then, in local coordinates, we can write the coordinates $(x_1, \dots, x_n)$ by
\begin{align*}
(x^{(0)},x^{(1)},\dots, x^{(\ell)}),
\end{align*}
where $x^{(j)}$ denotes the coordinates $(x_i)_{i\in \hat{I}_j}\ (j=0,\dots, \ell)$. 

Under the conditions H1-H3, we may obtain multi-normal deformation along $\chi$(see \cite[Section 1.2]{HPY} for details).
%Moreover we may multi-microlocalize subanalytic sheaves along the family of submanifolds $\chi$. 
 As usual we denote by $\tau_{M_1}: T_{M_1} X\to M_1$ the normal bundle. We first consider a normal deformation $\widetilde{X}_{M_1}$ along $M_1$(see \cite[Section 4.1]{KS90}), i.e. an analytic manifold $\widetilde{X}_{M_1}$, a pair of mappings $p_{M_1}: \widetilde{X}_{M_1} \to X, t_1: \widetilde{X}_{M_1} \to \mathbb{R}$, and an action of $\mathbb{R}\setminus \{0\}$ on $\widetilde{X}_{M_1}$ $(\widetilde{x}, r) \mapsto \widetilde{x}\cdot r$ such that $p^{-1}_{M_1}(X \setminus {M_1}) \simeq X \setminus M_1\times \mathbb{R}\setminus \{0\}$, $t_1^{-1}(c) \simeq X$ for each $c \not = 0$ and $t_1^{-1}(0) \simeq T_{M_1}X$.
We denote by $\Omega_{M_1}$ the open subset of $\widetilde{X}_{M_1}$ obtained as the inverse image of $\mathbb{R}^+$ by $t_1$, by $i_{\Omega_{M_1}}$ the embedding $\Omega_{M_1} \to \widetilde{X}_{M_1}$, $\tilde{p}_{M_1}= p_{M_1}\circ i_{\Omega_{M_1}}$ and $s_{M_1}: t_1^{-1}(0)\simeq T_{M_1} X \hookrightarrow \widetilde{X}_{M_1}$ the inclusion. We get the commutative diagram
% https://q.uiver.app/#q=WzAsNSxbMCwwLCJUX3tNXzF9WCJdLFsxLDAsIlxcd2lkZXRpbGRle1h9X3tNXzF9Il0sWzIsMCwiXFxPbWVnYV97TV8xfSJdLFsxLDEsIlguIl0sWzAsMSwiTV8xIl0sWzAsMSwic197TV8xfSJdLFsyLDEsImlfe1xcT21lZ2Ffe01fMX19IiwyXSxbMSwzLCJwX3tNXzF9Il0sWzIsMywiXFx0aWxkZXtwfV97TV8xfSJdLFswLDQsIlxcdGF1X3tNXzF9Il0sWzQsMywiaV97TV8xfSJdXQ==
\[
\begin{tikzcd}
	T_{M_1}X \arrow[r, "{s_{M_1}}"] \arrow[d, "{\tau_{M_1}}"'] & \widetilde{X}_{M_1} \arrow[r, "{i_{\Omega_{M_1}}}"] \arrow[d, "{p_{M_1}}"] & \Omega_{M_1} \arrow[dl, "{\tilde{p}_{M_1}}"] \\
	M_1 \arrow[r, "{i_{M_1}}"] & X
\end{tikzcd}
\]
Set $\widetilde{\Omega}_{M_1}= \{(x_1,t_1); t_1 \not=0\}$ and $\widetilde{M_2}:= \overline{(p_{M_1}|_{\widetilde{\Omega}_{M_1}})^{-1} M_2}$.
Then $\widetilde{M_2}$ is a closed smooth submanifold of $\widetilde{X}_{M_1}$.
Now we can define the normal deformation along $M_1$, $M_2$ as $\widetilde{X}_{M_1,M_2} := (\widetilde{X}_{M_1})_{\widetilde{M_2}}^{\sim}$.
Then we can define recursively the normal deformation along $\chi$ as $(\widetilde{X}_{M_1, \dots, M_{\ell-1}})_{\widetilde{M_\ell}}^{\sim}$.
Set $S_\chi = \{t_1, \dots, t_\ell =0\}$, $M = \bigcap_{i=1}^\ell M_i$ and $\Omega_\chi = \{t_1, \dots, t_\ell >0\}$. Then we have the commutative diagram
% https://q.uiver.app/?q=WzAsNSxbMCwwLCJTX1xcY2hpIl0sWzEsMCwiXFx3aWRldGlsZGV7WH0iXSxbMiwwLCJcXE9tZWdhX1xcY2hpIl0sWzEsMSwiWCJdLFswLDEsIk0iXSxbMCwxLCJzIiwyXSxbMiwxLCJpX1xcT21lZ2EiLDJdLFsxLDMsInAiXSxbMiwzLCJcXHRpbGRle3B9Il0sWzAsNCwiXFx0YXUiXSxbNCwzLCJpX00iXV0=
\[
\begin{tikzcd}
	S_\chi \arrow[r, "s"] \arrow[d, "\tau"'] & \widetilde{X} \arrow[r, "i_\Omega"] \arrow[d, "p"] & \Omega_\chi \arrow[dl, "\tilde{p}"] \\
	M \arrow[r, "i_M"] & X.
\end{tikzcd}
\]

Let us consider the canonical map $T_{M_j}\iota(M_j) \to M_j \hookrightarrow X, (j=1,\cdots, \ell)$ we write for short
\begin{equation*}S_\chi \simeq \underset{X, 1\leq j \leq \ell}{\times}T_{M_j}\iota(M_j):= T_{M_1}\iota(M_1) \underset{X}{\times}T_{M_2}\iota(M_2) \underset{X}{\times}\dots \underset{X}{\times}T_{M_\ell}\iota(M_\ell). 
\end{equation*}

\begin{defi}
The multi-specialization along $\chi$ is the functor
\[\function{\nu_\chi^{sa}}{D^b(k_{X_{sa}})}{D^b(k_{{S_\chi}_{sa}}),}{F}{s^{-1}R\Gamma_{\Omega_\chi} p^{-1}F.}\]
We also define the functor $\nu_\chi = \rho^{-1} \nu_\chi^{sa}: D^b(k_{X_{sa}}) \to D^b(k_{S_\chi})$.
\end{defi}
Let  $S^\ast_\chi$ be the dual vector bundle of $S_\chi$:
 \[S^\ast_\chi := T^\ast_{M_{1}}\iota(M_{1}) \times_X \cdots \times_X T^\ast_{M_{\ell}} \iota(M_{\ell}).\]
\begin{defi}
The multi-microlocalization along $\chi$ is the functor
\[\function{\mu_\chi^{sa}}{D^b(k_{X_{sa}})}{D^b(k_{S^\ast_{\chi}}\kern-1pt{}_{{}_{sa}}),}{F}{\nu_\chi^{sa}(F)^{\wedge_{S_\chi}}.}\]
We set $\mu_\chi := \rho^{-1} \mu_\chi^{sa}: D^b(k_{X_{sa}}) \to D^b(k_{S^\ast_\chi}).$
\end{defi}

Now we are going to restate a stalk formula for multi-microlocalization proved in \cite[Theorem 2.11]{HPY}. As the problem is local, we may assume that $X=\mathbb{R}^n$ and $q=0$ with coordinates $(x_1, \dots, x_n)$, and that there exists a subset $I_k(k=1,2,\dots, \ell)$ in $\{1,2,\dots, n\}$ with the conditions ((\ref{1.3}) (i),(ii)) such that each submanifold $M_k$ is given by $\{x=(x_1, \dots, x_n) \in \mathbb{R}^n; x_i =0(i\in I_k)\}$. Recall that $\hat{I}_k$ was defined by (\ref{ihat}) and that we set $M= \cap_k M_k$ and $n_k= \sharp\hat{I}_k$.
Then locally we have
\begin{align}\label{local}
X= M \times (N_1 \times N_2 \times \dots \times N_\ell)= M \times N,
\end{align}
where $N_k$ is $\mathbb{R}^{n_k}$ with coordinates $x^{(k)} = (x_i)_{i\in \hat{I}_k}$. Set, for $k \in \{1,\dots, \ell\}$,
\begin{align}
J_{\prec k}&:= \{j\in \{1, \dots,\ell\}, I_j \subsetneq I_k\}, \nonumber \\
J_{\succ k}&:= \{j\in \{1, \dots,\ell\},I_j \supsetneq I_k\},\\
J_{\nparallel k}&:= \{j \in \{1,\dots, \ell\}, I_j \cap I_k = \emptyset\}\nonumber.
\end{align}
Let $p=p_1 \times \cdots \times p_\ell = (q;\xi^{(1)}, \dots, \xi^{(\ell)}) \in T^\ast_{M_1}\iota(M_1) \times_X \dots \times_X T^\ast_{M_\ell}\iota(M_\ell)$ and consider the following conic subset in $N$
\begin{align}\label{gamma}
\gamma_k:=\left\{ (x^{(1)}, x^{(2)},\cdots,x^{(\ell)}) \in N\left|
\begin{array}{ll}
x^{(j)}=0,\  ( j \in J_{\prec k} \sqcup J_{\nparallel k}), \\
x^{(j)}\in \mathbb{R}^{n_j},\  (j\in J_{\succ k}),\\
\langle x^{(j)} ,\xi^{(k)} \rangle > 0,\  (j=k)
\end{array}
\right. \right\}.
\end{align}
\begin{theo}\label{stalkmu}
\textit{
Let $p = p_1 \times \cdots \times p_\ell = (q; \xi^{(1)}, \dots, \xi^{(\ell)}) \in S^\ast_\chi$, and let $F \in D^b(k_{X_{sa}})$. Then we have
\begin{align*}
H^k(\mu_\chi F)_p \simeq \varinjlim_{G, U} H^k_G(U; F).
\end{align*}
Here $U$ is an open subanalytic neighborhood of $q$ in $X$ and $G$ is a closed subanalytic subset in the form $M \times \sum_{k=1}^\ell G_k$ with $G_k$ being a closed subanalytic convex cone in $N$ satisfying $G_k \setminus \{0\}\subset \gamma_k$, where $\gamma_k$ is defined in (\ref{gamma}).
}
\end{theo}

In \cite[Theorem 4.1]{HPY}, the authors proved the existence of Uchida's triangle(\cite[Appendix A]{U}) in the context of multi-microlocal analysis. The theorem only applies to classical sheaves, but we have shown that it holds for sheaves on subanalytic sites as well.

\begin{theo}\label{ut}
{\it
Let $F \in D^b_{(\mathbb{R}^+)^\ell}(k_{E_{sa}})$. There exists a natural isomorphism
\[\tau^!R\tau_{!!} F\simeq Rp_{1\ast}p^{!}_2(F^{\wedge_E}),\]
and the natural morphism $F \to \tau^! R\tau_{!!} F$ is embedded to the following distinguished triangle:
% https://q.uiver.app/?q=WzAsNCxbMCwwLCJGIl0sWzEsMCwiXFx0YXVeISBSXFx0YXVfISBGIl0sWzIsMCwiUnBfezFcXGFzdH1eKyBwXzJeeyshfShGXntcXHdlZGdlX0V9KSJdLFszLDAsInt9Il0sWzAsMV0sWzEsMl0sWzIsMywiKzEiXV0=
\[
\begin{tikzcd}
	F \arrow[r] & \tau^! R\tau_{!!} F \arrow[r] & Rp_{1\ast} R\Gamma_{P^+}(p_2^{!}(F^{\wedge_E})) \arrow[r, "{+1}"] & {}.
\end{tikzcd}
\]
}
\end{theo}

Next we explain the notion of the $\mu$hom in multi-microlocalizations(see \cite{H}).
Let $\chi$ satisfy the conditions H1-H3. Let $A$ be an $\ell\times \ell$ matrix by
\[a_{ji}=
\left\{
\begin{array}{ll}
1, & (i \in I_j )\\
0, & (i \not \in I_j).
\end{array}
\right.
\]
We call $A$ an associated matrix with $\chi$.
$X^2 = X \times X$ and $\Delta_X \subset X^2$ denotes its diagonal, i.e., $\Delta_X = \{(x,y)\in X^2; x=y\}$.

Define the ambient space $\widehat{X}$ of the universal family as
\[(X^2)^\ell= X^2 \times \cdots \times X^2.\]
Define the universal family $\widehat{\chi}= \{\widehat{M}_1, \dots, \widehat{M}_\ell\}$ of closed submanifolds in $\widehat{X}$ in the following manner:
Each $\widehat{M}_j$ has a form
\[\widehat{M}_j = Z_{j1} \times Z_{j2} \times \cdots \times Z_{j\ell} \subset \widehat{X},\]
where $Z_{jk} \subset X^2$ is either $X^2$ or $\Delta_X$ determined by the rules:
\[Z_{jk} := \left\{
\begin{array}{ll}
X^2, &  (a_{jk} =0), \\
\Delta_X, &  (a_{jk} = 1).
\end{array}
\right.
\]
For the universal family $\widehat{\chi}$, we may consider the multi-normal deformation $\widetilde{(X^2)^\ell}_{\widehat{\chi}}$ of $(X^2)^\ell$ along $\widehat{\chi}$ and the multi-microlocalization $\mu^{sa}_{\widehat{\chi}}$(see \cite{H} for details.). And we have
\[
S_{\widehat{\chi}}^\ast = T^\ast X \times \cdots \times T^\ast X.
\]
\begin{defi}[Multi-microlocal Hom functors]
Let $F_1, \dots, F_\ell$ and $G$ be in $D^b(X_{sa})$. We define the functor as 
%from $(D^b(X_{sa})^\circ)^\ell \times D^b(X_{sa})$ to $D^+(k_{S^\ast_{\chi}}\kern-1pt{}_{{}_{sa}})$ as
\[\mu hom^{sa}_{\widehat{\chi}}(F_1, \dots, F_\ell; G):= \mu^{sa}_{\widehat{\chi}}(\RintHom(p^{-1}_2 (F_1 \boxtimes \cdots \boxtimes F_\ell), p^!_1 Ri_{\Delta\ast} G))\]
and we set the functor 
%$(D^b(X_{sa})^\circ)^\ell \times D^b(X_{sa})$ to $D^b(k_{S^\ast_{\chi}})$ 
as
\[\mu hom_{\widehat{\chi}}(F_1, \dots, F_\ell; G):= \rho^{-1}\mu hom^{sa}_{\widehat{\chi}}(F_1, \dots, F_\ell; G),\]
where
% https://q.uiver.app/?q=WzAsNCxbMCwwLCJYXmwiXSxbMSwwLCJcXHdpZGVoYXR7WH0iXSxbMiwwLCJYXmwiXSxbMywwLCJYIl0sWzEsMCwicF8xIiwyXSxbMSwyLCJwXzIiXSxbMywyLCJpX1xcRGVsdGEiLDJdXQ==
\[
\begin{tikzcd}
	X^\ell & \widehat{X} \arrow[l, "{p_1}"'] \arrow[r, "{p_2}"] & X^\ell & X \arrow[l, "{i_\Delta}"']
\end{tikzcd}
\]
\begin{enumerate}
\item[$-$] $i_\Delta: X \hookrightarrow X^\ell$ is defined by $x \in X \to (x,x, \dots, x) \in X^\ell$,
\item[$-$] $p_1: (X^2)^\ell \to X^\ell$ is the canonical projection $(x_1, y_1) \times \cdots \times (x_\ell, y_\ell)\to (x_1, \dots, x_\ell)$ and
\item[$-$] $p_2: (X^2)^\ell \to X^\ell$ is the canonical projection $(x_1, y_1) \times \cdots \times (x_\ell, y_\ell)\to (y_1, \dots, y_\ell)$.
\end{enumerate}
\end{defi}
\section{A stalk formula for the multi-microlocal $\Hom$ functors}\label{sec:stalkformula}
In this sectin we are going to prove a stalk formula for the multi-microlocal Hom functors. As the question being local, we may assume $X$ is a vector space $\mathbb{R}^n$.
\begin{lem}\label{gamlem}
{\it
Let $F \in D^b_{\mathbb{R}\text{-}c}(k_{X})$ and $V$ be a relatively compact open subanalytic set. For a closed convex cone $\gamma$ in $X$, we have 
\[Rp_{1!!}(p_2^{-1}(R\rho_\ast F)_V)_{Z(\gamma)}\simeq R\rho_\ast ( \phi^{-1}_{\tilde{\gamma}}R\phi_{\gamma \ast}(F_V)),\]
where $p_1$ and $p_2$ the first and second projection from $X \times X$ to $X$.
}
\end{lem}
\proof
This follows immediately from \cite[Proposition 2.2.1]{Lu1}, \cite[Proposition 1.3.3]{Lu1}, \cite[Proposition 1.3.2]{Lu1} and \cite[Proposition 3.5.4]{KS90}.
In fact,
\begin{align}
Rp_{1!!}(p_2^{-1}(R\rho_\ast F)_V)_{Z({\gamma})} \simeq& Rp_{1!!}(R\rho_\ast(p_2^{-1}(F_V)_{Z({\gamma})})\label{iso1}\\
\simeq& R\rho_\ast Rp_{1!}(p_2^{-1}(F_V))_{Z({\gamma})}\label{2}\\
\simeq& R\rho_\ast (\phi^{-1}_{{\gamma}}R\phi_{{\gamma}\ast}(F_V))\nonumber.
\end{align}
The isomorphism (\ref{iso1}) follows since $F$ is constructible. The isomorphism (\ref{2}) follows since $p_1$ is proper on $\supp (p_2^{-1}(F_V))_{Z({\gamma})}$.
\endproof
\begin{theo}\label{stalkformula}
\textit{
Let $p=p_1 \times p_2 \times \cdots \times p_\ell=(q; \xi^{(1)}, \dots, \xi^{(\ell)}) \in \overbrace{T^*X \times_X \cdots \times_X T^\ast X}^{\ell}$. 
$G_1, \dots, G_\ell \in \mbox{Ob}(D^b_{\mathbb{R}\text{-}c}(k_{X}))$ and $F \in \mbox{Ob}(D^b(k_{X_{sa}}))$ we have
\begin{align*}{}&H^j\left(\mu hom_{\widehat{\chi}} (R\rho_\ast G_1, R\rho_\ast G_2, \dots, R\rho_\ast G_\ell; F)\right)_p \\
\simeq& \varinjlim_{U, \tilde{\gamma}} H^j\left(R\Gamma\left(U; \RintHom(R\rho_\ast (i^{-1}_{\Delta}\phi_{\tilde{\gamma}}^{-1}R\phi_{\tilde{\gamma}\ast}((G_1\boxtimes\cdots \boxtimes G_\ell)_{U \times \cdots \times U})),F)\right)\right),
\end{align*}
where
$\tilde{\gamma}=\tilde{\gamma_1} + \cdots + \tilde{\gamma_\ell}\subset X^\ell$, each $\tilde{\gamma_k}$ is a closed conic subanalytic proper convex subset of $\gamma_k\cup \{0\}$,
\[
\gamma_k:=\left\{(y^{(1)}, \cdots, y^{(\ell)}) \in X^\ell \left|
\begin{array}{ll}
y^{(j)}=0, ( j \in J_{\prec k} \sqcup J_{\nparallel k}), \\
y^{(j)}\in \mathbb{R}^{n}, (j\in J_{\succ k}),\\
\langle y^{(j)}, \xi^{(k)} \rangle < 0, (j=k)
\end{array}
\right. \right\}
\]
and $U$ ranges through the family of open subanalytic neighborhoods of $q$.
}
\end{theo}
\proof
 As the question being local, we may assume $X$ is a vector space $\mathbb{R}^n$ and $q =0$ with coordinates $(x_1, \dots, x_n)$. Clearly (\ref{local}) becomes $X^{2\ell} = X^\ell \times X^\ell$.
Set
\[
\delta_k:=\left\{ (y^{(1)}, \cdots, y^{(\ell)}) \in X^\ell \left|
\begin{array}{ll}
y^{(j)}=0, ( j \in J_{\prec k} \sqcup J_{\nparallel k}), \\
y^{(j)}\in \mathbb{R}^{n}, (j\in J_{\succ k}),\\
\langle y^{(j)}, \xi^{(k)} \rangle > 0, (j=k)
\end{array}
\right. \right\}
.\]
We propose the following identification $\varphi$.
\begin{align}\label{id}
\idfunction{X^\ell \times X^\ell}{X^\ell \times X^\ell}{(x_1,x_2, \dots, x_\ell ;y_1, y_2, \dots, y_\ell)}{\left(x_1 ,x_2, \dots, x_\ell ; x_1- y_1, x_2- y_2,\dots, x_\ell-y_\ell\right)}{\varphi}
\end{align}
We have the following claim.
\begin{claim}
Let $G$ be a closed subanalytic subset in the form $ X^\ell \times \sum_{k=1}^\ell G_k$ with $G_k$ being a closed subanalytic proper convex cone in $X^\ell$ satisfying $G_k\setminus \{0\} \subset \delta_k$. Then for each such $G$ there exists a closed conic subanalytic proper convex  subset $\tilde{\gamma_k} \setminus \{0\} \subset -\delta_k = \gamma_k$ such that
\[ G \subset  \varphi^{-1}\left(Z({\tilde{\gamma}})\right),\]
where $\tilde{\gamma}= \tilde{\gamma_1} +\tilde{\gamma_2}+\dots +\tilde{\gamma_\ell}$ and
\[
Z({\tilde{\gamma}}): = \left\{(x_1, x_2, \dots, x_\ell; y_1,y_2,\dots, y_\ell)\in X^\ell \times X^\ell\left|\right. (y_1 -x_1,y_2-x_2, \dots, y_\ell-x_\ell) \in \tilde{\gamma}\right\}.
\]
Moreover for such $\tilde{\gamma}$ we may find a closed subanalytic subset $G^\prime= X^\ell \times \sum_{k=1}^\ell G^\prime_k$ with $G^\prime_k$ being a closed subanalytic proper convex  cone in $X^\ell$ satisfying $G^\prime_k\setminus \{0\} \subset \delta_k$ so that $\varphi^{-1}\left(Z({\tilde{\gamma}})\right) \subset G^\prime$.
\end{claim}
In fact, let $G$ be such a conic set. Then there exists $G_k \setminus \{0\} \subset \delta_k$. We may find $\tilde{\gamma_k}$ such that $G_k \subset -\tilde{\gamma_k}$ and $-\tilde{\gamma_k}\setminus \{0\} \subset -\delta_k$. By the identification (\ref{id}), $\varphi^{-1}(Z({\tilde{\gamma}})) = \left\{(x; y )\in X^\ell \times X^\ell \left|\right. -y \in {\tilde{\gamma}}\right\}$. If $(x, y) \in G$ then $y \in \sum_{k=1}^\ell G_k $. And $y \in -\tilde{\gamma}$ thus $(x,y)\in \varphi^{-1}(Z({\tilde{\gamma}}))$. Finally take a closed subanalytic proper convex cone $G_k^\prime$ in $X^\ell$ so that $-\tilde{\gamma_k}\subset G^\prime_k$ and $G^\prime_k\setminus \{0\} \subset \delta_k$ for each $k$. It is easy to see  $\varphi^{-1}(Z({\tilde{\gamma}})) \subset G^\prime= X^\ell \times \sum_{k=1}^\ell G^\prime_k$.
Therefore by Theorem \ref{stalkmu} and Lemma \ref{gamlem} we obtain the chain of isomorphisms (put $\overset{\ell}{\underset{i=1}{\boxtimes} }G_i:= G_1\boxtimes \cdots \boxtimes G_\ell$)
\begin{align*}
 {}&H^k \mu hom_{\widehat{\chi}}(R\rho_\ast G_1, R\rho_\ast G_2, \dots, R\rho_\ast G_\ell; F)_p \\
\simeq& \ H^k \mu_{\widehat{\chi}} (\RintHom(p_2^{-1}(R\rho_\ast \overset{\ell}{\underset{i=1}{\boxtimes} }G_i), p^!_1Ri_{\Delta\ast}F ))_{p}\\
\simeq& \varinjlim_{G,U,V}H^k(R\Gamma_G (U\times V; \RintHom(p_2^{-1}(R\rho_\ast \overset{\ell}{\underset{i=1}{\boxtimes} }G_i), p^!_1Ri_{\Delta\ast}F)))\\
\simeq& \varinjlim_{{\tilde{\gamma}},U,V}H^k(R\Gamma_{Z({\tilde{\gamma}})} (U\times V; \RintHom(p_2^{-1}(R\rho_\ast \overset{\ell}{\underset{i=1}{\boxtimes} }G_i), p^!_1Ri_{\Delta\ast}F))\\
\simeq& \varinjlim_{\tilde{\gamma},U,V}H^k(R\Gamma(U\times X; \RintHom((p_2^{-1}(R\rho_\ast \overset{\ell}{\underset{i=1}{\boxtimes} }G_i)_V)_{Z({\tilde{\gamma}})} , p^!_1Ri_{\Delta\ast}F)))\\
\simeq& \varinjlim_{\tilde{\gamma},U,V}H^k(R\Gamma(U; \RintHom(Rp_{1!!}((p_2^{-1}(R\rho_\ast \overset{\ell}{\underset{i=1}{\boxtimes} }G_i)_V)_{Z({\tilde{\gamma}})}) , Ri_{\Delta\ast}F)))\\
\simeq& \varinjlim_{\tilde{\gamma},U, V}H^k(R\Gamma(i^{-1}_\Delta U; \RintHom(i^{-1}_\Delta Rp_{1!!}(p_2^{-1}(R\rho_\ast \overset{\ell}{\underset{i=1}{\boxtimes} }G_i )_V)_{Z({\tilde{\gamma}})} ,  F)))\\
\simeq& \varinjlim_{\tilde{\gamma},U}H^k(R\Gamma(i^{-1}_\Delta U; \RintHom(i^{-1}_\Delta R\rho_\ast\phi^{-1}_{\tilde{\gamma}}R\phi_{\tilde{\gamma}\ast}((\overset{\ell}{\underset{i=1}{\boxtimes} }G_i)_U) ,  F))).
\end{align*}
$G$ runs through a family of closed subanalytic subsets in the form $ X^\ell \times \sum_{k=1}^\ell G_k$ with $G_k$ being a closed subanalytic proper convex cone in $X^\ell$ satisfying $G_k\setminus \{0\} \subset \delta_k$, $U$ and $V$ runs through a family of open subanalytic neighborhoods of $i_\Delta(q)$ in $X^\ell$. Notice that $U$ is different from the one in the statement of theorem. The first isomorphism follows from \cite[Proposition 1.3.3]{Lu1} and \cite[Proposition 1.3.1]{Lu2}.
\endproof
\begin{lem}\label{lemgm}
\textit{
Let $\gamma$ be a closed proper convex cone (with vertex at 0) in a vector space $V$, $G$ be a closed convex subset. 
If $\gamma^a + G$ is closed then
\[
\phi^{-1}_{\gamma}R\phi_{\gamma\ast} \mathbb{C}_{G} \simeq \mathbb{C}_{\gamma^a + G}.
\]
Here we denote by ``$a$'' the antipodal map $x \mapsto -x$.
}
\end{lem}
\proof
Let us first prove $\phi^{-1}_{\gamma}R\phi_{\gamma\ast} \mathbb{C}_{G+\gamma^a} \xrightarrow{\sim} \mathbb{C}_{G+\gamma^a}$.
By adjunction, there exists a morphism $\phi^{-1}_{\gamma}R\phi_{\gamma\ast}\mathbb{C}_{G+\gamma^a} \to \mathbb{C}_{G+\gamma^a}$. Then taking a stalk of the cohomology proves the assertion. In fact, by \cite[Theorem 3.5.3]{KS90} for any $x \in X$
\[
\varinjlim_{x\in U} H^k(U;\phi^{-1}_{\gamma}R\phi_{\gamma\ast} \mathbb{C}_{{G+\gamma^a}}) \simeq \varinjlim_{x\in U} H^k(U+\gamma;\mathbb{C}_{G+\gamma^a})
\]
and since $(x+ \gamma) \cap G+\gamma^a$ is convex when $(x + \gamma) \cap G+\gamma^a \ne \emptyset$ we have $R\Gamma((x+\gamma) \cap G+\gamma^a; \mathbb{C}_{G+\gamma^a}) \simeq \mathbb{C}$.
Thus it is enough to prove if $x \not \in G+\gamma^a$ then $(x+\gamma) \cap G+\gamma^a = \emptyset$. By contradiction, if we may find $g_1, g_2 \in \gamma$ and $ g_3 \in G$ so that $x+ g_1 = g_3 -g_2$ then $x =g_3 -g_1-g_2$. This implies $x \in G+\gamma^a$ which is a contradiction.

Now we prove the general assertion. We have a morphism $\mathbb{C}_{G + \gamma^a} \to \mathbb{C}_{G}$. Therefore by the first assertion we get a morphism
\[
\mathbb{C}_{G + \gamma^a} \simeq \phi^{-1}_{\gamma}R\phi_{\gamma\ast}  \mathbb{C}_{G + \gamma^a} \to \phi^{-1}_{\gamma}R\phi_{\gamma\ast}  \mathbb{C}_{G}.
\]
Then taking a stalk of the cohomology proves the assertion. In fact, by \cite[Theorem 3.5.3]{KS90}
\[
\varinjlim_{x\in U} H^k(U;\phi^{-1}_{\gamma}R\phi_{\gamma\ast} \mathbb{C}_{{G}}) \simeq \varinjlim_{x\in U} H^k(U+\gamma;\mathbb{C}_{G})
\]
and since $(x+ \gamma) \cap G $ is convex when $(x + \gamma) \cap G  \ne \emptyset$ we have $R\Gamma((x+\gamma) \cap G ; \mathbb{C}_{G}) \simeq \mathbb{C}$.
Thus it is enough to prove if $x \not \in G +\gamma^a$ then $(x+\gamma) \cap G  = \emptyset$. By contradiction, if we may find $g_1 \in \gamma$ and $ g_2 \in G$ so that $x+ g_1 = g_2 $ then $x =g_2 -g_1$. This implies $x \in G+\gamma^a$ which is a contradiction.
\endproof
\begin{cor}[Normal crossing type]
\textit{
Assume $X$ is a vector space.
Let $\widehat{\chi}=\{\widehat{D}_1, \dots, \widehat{D}_\ell\}$ where
\begin{align*}
\widehat{D}_j =& Z_{j1} \times Z_{j2} \times \cdots \times Z_{j\ell} \subset \widehat{X},\\
Z_{jk} :=& \left\{
\begin{array}{ll}
\Delta_X, &  (j=k), \\
X^2, &  (o.w.).
\end{array}
\right.
\end{align*}
Let $M_1,\dots, M_\ell$ be $\ell$ closed subanalytic convex subsets of $X$.   Then for $p =(q; \xi_1, \dots, \xi_\ell) \in T^\ast X \times_X T^\ast X \times_X \cdots \times_X T^\ast X$ with $M_k \subset \{x \in X; \langle x-q, \xi_k \rangle \ge 0\}(k=1, \dots, \ell)$, we have
\begin{align*}
{}&H^j(\mu hom_{\widehat{\chi}} (\mathbb{C}_{M_1}, \mathbb{C}_{M_2}, \dots, \mathbb{C}_{M_\ell}; F))_p \\
\simeq& \varinjlim_{U, \widetilde{\delta_1},\cdots, \widetilde{\delta_\ell} } H^j_{(M_1 + \widetilde{\delta_1} )\cap (M_2 + \widetilde{\delta_2}) \cap \cdots \cap (M_\ell + \widetilde{\delta_\ell})\cap U}(U; F),
\end{align*}
where
each $\widetilde{\delta_k}$ is a closed subanalytic proper convex cone of $\delta_k\cup \{0\}$,
$\delta_k =\{x \in X; \langle x, \xi_k \rangle > 0\}$.
}
\end{cor}
\proof
By Theorem \ref{stalkformula} we have 
\begin{align*}
{}& H^j(\mu hom_{\widehat{\chi}} (\mathbb{C}_{M_1}, \dots, \mathbb{C}_{M_\ell}; F))_p \\
\simeq& \varinjlim_{U, \tilde{\gamma}} H^j(U; \RintHom(R\rho_\ast i^{-1}_{\Delta}(\phi_{\tilde{\gamma}}^{-1}R\phi_{\tilde{\gamma}\ast}(\mathbb{C}_{M_1\cap U\times M_2\cap U \times\cdots \times M_\ell\cap U})),F))),
\end{align*}
where
$\tilde{\gamma}=\tilde{\gamma}_1 + \cdots + \tilde{\gamma}_\ell\subset X^\ell$, each $\tilde{\gamma_k}$ is a closed conic subanalytic proper convex subset of $\gamma_k\cup \{0\}$,
\[
\gamma_k:=\left\{(x^{(1)}, \cdots, x^{(\ell)}) \in X^\ell \left|
\begin{array}{ll}
\langle x^{(j)}, \xi_k \rangle < 0, (j=k), \\
x^{(j)}=0, (\mbox{o.w.})
\end{array}
\right. \right\}
\]
and $U$ ranges through the family of open subanalytic neighborhoods of $q$ in $X$. 

Let $\tilde{\gamma}$ be a cone as above. We may find an open neighborhood $U$ of $q$ such that:
$$M_1\cap U\times M_2\cap U \times \cdots \times M_\ell \cap U = M_1\times M_2 \times \cdots \times M_\ell \cap (U \times \cdots \times U + \tilde{\gamma}).$$
Thus for such $U$ we have:
\[
(\phi_{\tilde{\gamma}}^{-1}R\phi_{\tilde{\gamma}\ast} \mathbb{C}_{M_1\cap U\times M_2\cap U \times \cdots \times M_\ell \cap U})|_{U \times \cdots \times U } \simeq (\phi_{\tilde{\gamma}}^{-1}R\phi_{\tilde{\gamma}\ast} \mathbb{C}_{M_1\times M_2 \times \cdots \times M_\ell} )|_{U \times \cdots \times U }.
\]
Therefore by Lemma \ref{lemgm} we have, 
\begin{align*}
{}& H^j(\mu hom_{\widehat{\chi}} (\mathbb{C}_{M_1}, \dots, \mathbb{C}_{M_\ell}; F))_p \\
\simeq& \varinjlim_{U, \tilde{\gamma}} H^j(U; \RintHom(R\rho_\ast i^{-1}_{\Delta}(\phi_{\tilde{\gamma}}^{-1}R\phi_{\tilde{\gamma}\ast}(\mathbb{C}_{M_1\cap U\times M_2\cap U \times \cdots \times M_\ell \cap U})),F)))\\
\simeq& \varinjlim_{U, \tilde{\gamma}} H^j(U; \RintHom(R\rho_\ast i^{-1}_{\Delta}(\mathbb{C}_{M_1\times M_2 \times \cdots \times M_\ell+ \tilde{\gamma}\cap {U \times \cdots \times U}}),F))\\
    \simeq& \varinjlim_{U, \tilde{\gamma}} H^j(U; \RintHom(R\rho_\ast (i^{-1}_{\Delta}(\mathbb{C}_{(M_1 + q_1 (\tilde{\gamma_1}^a ))\times (M_2  + q_2 (\tilde{\gamma_2}^a))\times \cdots\times (M_\ell  + q_\ell (\tilde{\gamma_\ell}^a))\cap {U \times \cdots \times U}})),F))\\
\simeq& \varinjlim_{U, \tilde{\gamma}} H^j(U; \RintHom(R\rho_\ast \mathbb{C}_{(M_1 + q_1 (\tilde{\gamma_1}^a ))\cap (M_2 + q_2 (\tilde{\gamma_2}^a))\cap \cdots\cap (M_\ell+ q_\ell (\tilde{\gamma_\ell}^a))\cap U},F))\\
\simeq& \varinjlim_{U, \tilde{\gamma} } H^j_{(M_1 + q_1 (\tilde{\gamma_1}^a ))\cap (M_2 + q_2 (\tilde{\gamma_2}^a))\cap \cdots\cap (M_\ell+ q_\ell (\tilde{\gamma_\ell}^a))\cap U}(U; F),
\end{align*}
where $q_k$ is $k$-th projection. 
\endproof
\begin{ex}
Let $X=\mathbb{C}^2$ and $\widehat{\chi}= \{\Delta \times \mathbb{C}^2 , \mathbb{C}^2 \times \Delta \}$. Consider closed subsets
$M_1 = \{ y_1 = 0, x_1 \ge 0\} = \mathbb{R}_{\ge 0} \times  \mathbb{C}$ and
$M_2 = \{ y_2 = 0, x_2 \ge 0\} =  \mathbb{C}   \times \mathbb{R}_{\ge 0}$.

%$\mathscr{O}^t_X$
By Corollary \ref{stalkformula} for any subanalytic sheaf $F \in D^b(\mathbb{C}_{X_{sa}})$ and $p = (q; \zeta_1, \zeta_2) \in T^\ast X \times_X T^\ast X = T^\ast \mathbb{C}^2 \times_{\mathbb{C}^2} T^\ast \mathbb{C}^2 =\overbrace{\underbrace{ \mathbb{C}^2}_{\text{base space}}\times \underbrace{\mathbb{C}^2 \times \mathbb{C}^2}_{\text{fiber space}}}^{\text{total space}}$ with $M_k \subset \{z \in X; \mbox{Re}\langle z-q, \zeta_k \rangle \ge 0\}(k=1,2)$,

\begin{align*}
{}&H^j(\mu hom_{\widehat{\chi}} (\mathbb{C}_{M_1}, \mathbb{C}_{M_2}; F))_p \\
\simeq& \varinjlim_{U, \widetilde{\delta_1}, \widetilde{\delta_2} } H^j_{(M_1 + \widetilde{\delta_1} )\cap (M_2 + \widetilde{\delta_2})\cap U}(U; F),
\end{align*}
where
each $\tilde{\delta_k}$ is a closed subanalytic proper convex cone of $\delta_k\cup \{0\}$,
$\delta_k =\{z\in \mathbb{C}^2; \mbox{Re}\langle z, \zeta_k \rangle > 0\}$
and $U$ ranges through the family of open subanalytic neighborhoods of $q$. Here $\langle (z_1, z_2),(z_3, z_4) \rangle := z_1 z_3 +z_2z_4$.

Let $p = ((1,1); (\sqrt{-1}, 0), (0,\sqrt{-1}))$, we have
\begin{align*}
{}& H^j(\mu hom_{\widehat{\chi}} (\mathbb{C}_{M_1}, \mathbb{C}_{M_2}; F))_p \\
\simeq& \varinjlim_{U} H^j_{H \times H \cap U}(U; F),
\end{align*}
where
%$\mathbb{C}^2 \times \mathbb{C}^2 \ni (z_{11}, z_{12}, z_{21},z_{22})= (x_{11}+iy_{11}, x_{12}+iy_{12}, x_{21}+iy_{21}, x_{22}+iy_{22})$
$U$ ranges through the family of open subanalytic neighborhoods of $(1,1)$ and $H:=\{x+iy; y\leq 0\}$. If $F=R\rho_\ast \mathscr{O}_X$ then by \cite[Theorem 2.2.2]{KKK} or the useful criterion \cite[(2.9.14)]{KS90} for $j\not =2$ we have
\begin{align*}
{}& H^j(\mu hom_{\widehat{\chi}} (\mathbb{C}_{M_1}, \mathbb{C}_{M_2}; R\rho_*\mathscr{O}_X))_p \\
\simeq& \varinjlim_{U} H^j_{H \times H \cap U}(U; R\rho_\ast\mathscr{O}_X)= 0.
\end{align*}

When $p=((1, 0); (\sqrt{-1}, 0), (0,\sqrt{-1}))$
\begin{align*}
{}& H^j(\mu hom_{\widehat{\chi}} (\mathbb{C}_{M_1}, \mathbb{C}_{M_2}; F))_p \\
\simeq&  \varinjlim_{\lambda, U}H^j_{G_{\lambda}   \times H \cap U}(U, F),
\end{align*}
where
$G_\lambda = \lambda + \mathbb{R}_{\ge 0}\subset \mathbb{C}$ and $\lambda$ runs through the family of closed conic subanalytic proper convex subsets in $\mathbb{C}$ with its vertex being $0$ such that $\lambda \subset  \{y_1 < 0\} \cup \{0\}$ and $U$ ranges through the family of open subanalytic neighborhoods of $(1,0)$.  If $F=R\rho_*\mathscr{O}_X$ then by the useful criterion \cite[(2.9.14)]{KS90} for $j\not =2$ we have
\begin{align*}
 H^j(\mu hom_{\widehat{\chi}} (\mathbb{C}_{M_1}, \mathbb{C}_{M_2}; R\rho_*\mathscr{O}_X))_p = 0.
\end{align*}

Finally consider $p=((0, 0); (\sqrt{-1}, 0), (0,\sqrt{-1}))$.
\begin{align*}
{}& H^j(\mu hom_{\widehat{\chi}} (\mathbb{C}_{M_1}, \mathbb{C}_{M_2}; F))_p \\
\simeq&  \varinjlim_{\lambda_1, \lambda_2, U}H^j_{G_{\lambda_1}   \times G_{\lambda_2} \cap U}(U, F),
\end{align*}
where each $G_{\lambda_k} = \lambda_k + \mathbb{R}_{\ge 0}\subset \mathbb{C}$ and $\lambda_k$ runs through the family of closed conic subanalytic proper convex subsets in $\mathbb{C}$ with its vertex being $0$ such that $\lambda_k \subset  \{y_k < 0\} \cup \{0\}$ and $U$ ranges through the family of open subanalytic neighborhoods of $(0,0)$.  If $F=R\rho_*\mathscr{O}_X$ then by the useful criterion \cite[(2.9.14)]{KS90} for $j\not =2$ we have
\begin{align*}
 H^j(\mu hom_{\widehat{\chi}} (\mathbb{C}_{M_1}, \mathbb{C}_{M_2}; R\rho_*\mathscr{O}_X))_p = 0.
\end{align*}
\end{ex}

\section{Multi-microlocal Sato's triangle}\label{sec:satotr}
In this section we construct multi-microlocal Sato's triangle. Suppose $\chi= \{M_1, M_2, \dots, M_\ell\}$ satisfies the conditions H1, H2 and H3.
\begin{prop}\label{muiso}
\textit{
Let $i:M \to X$ be the canonical embedding and $\pi: S^\ast_\chi \to M$ be the canonical projection where $M=\bigcap_{j=1}^\ell M_j$. For $F \in D^b(k_{X_{sa}})$
\[\mu^{sa}_\chi(F)|_{M}\simeq R\pi_\ast \mu^{sa}_\chi(F) \simeq R\Gamma_{M} (F)|_{M}\simeq i^! F.\]
}
\end{prop}
\proof
It is enough to prove the second isomorphism. By Theorem \ref{ut}, for $H=\nu^{sa}_\chi(F)$ we have
\[\tau^! R\tau_{!!} H \simeq Rp_{1\ast}p^!_2 H^{\wedge_{S_\chi }}.\]
By the base change formula(\cite[Proposition 2.2.9]{Lu1}) on the Cartesian square:
% https://q.uiver.app/?q=WzAsNCxbMCwwLCJFIFxcdGltZXNfTSBFXlxcYXN0Il0sWzEsMCwiRV5cXGFzdCJdLFswLDEsIkUiXSxbMSwxLCJNIl0sWzAsMSwicF8yIl0sWzAsMiwicF8xIiwyXSxbMiwzLCJcXHRhdSIsMl0sWzEsMywiXFxwaSJdXQ==
\[
\begin{tikzcd}
	S_\chi \times_M S_\chi^\ast \arrow[r, "{p_2}"] \arrow[d, "{p_1}"'] & S_\chi^\ast \arrow[d, "\pi"] \\
	S_\chi \arrow[r, "\tau"'] & M
\end{tikzcd}
\]

we have $Rp_{1\ast}p^!_2 H^{\wedge_{S_\chi }} \simeq \tau^!R\pi_\ast H^{\wedge_{S_\chi }}$.
Thus
\[ \tau^! R\tau_{!!} H \simeq \tau^! R\pi_\ast H^{\wedge_{S_\chi }}.\]
Applying $i^!$ both sides:
\[R\tau_{!!} H \simeq R\pi_\ast H^{\wedge_{S_\chi }}.\]
Since $R\tau_{!!} \nu^{sa}_\chi(F) \simeq R\Gamma_M(F)$ by \cite[Proposition 2.6]{HPY}, we get
\[R\Gamma_M(F)\simeq R\pi_\ast \mu^{sa}_\chi(F).\]
\endproof

\begin{conv}
Since if we let $i: T^*X \times_X \cdots \times_X T^*X \to S^*_{\widehat{\chi}}$ be a canonical embedding we get by \cite[Lemma 2.17]{H}, $Ri_*i^{-1}\mu hom^{sa}(G_1, \dots, G_\ell; F) \simeq \mu hom^{sa}(G_1, \dots, G_\ell; F)$, $\mu hom^{sa}(G_1, \dots, G_\ell; F)$ may be seen as an object on $T^*X \times_X \cdots \times_X T^*X $. We take this convention throughout this paper.
\end{conv}
\begin{prop}
\textit{
Let $\pi: T^\ast X \times_X \cdots \times_X T^\ast X \to X$ be a canonical projection and $G_1, \dots, G_\ell, F \in \mbox{Ob}(D^b(k_{X_{sa}}))$.  Then 
\begin{align*}
R\pi_\ast \mu hom^{sa}_{\widehat{\chi}} (G_1,\dots, G_\ell; F) \simeq& \RintHom(G_1\boxtimes \cdots \boxtimes G_\ell, Ri_{\Delta\ast} F)\\
 \simeq& Ri_{\Delta\ast}\RintHom(G_1\otimes \cdots \otimes G_\ell, F).
\end{align*}
In particular,
\[\Hom_{D^b(k_{X_{sa}})}(\overset{\ell}{\underset{i=1}{\otimes} }G_i, F) \simeq H^0(T^\ast X \times_X \cdots \times_X T^\ast X; \mu hom^{sa}_{\widehat{\chi}} (G_1, \dots, G_\ell; F)).\]
}
\end{prop}
\proof
We have the following chain of isomorphisms:
\begin{align*}{}& R\pi_\ast \mu^{sa}_{\widehat{\chi}} \RintHom(q^{-1}_2 (G_1 \boxtimes \cdots \boxtimes G_\ell), q^!_1 Ri_{\Delta \ast}F)\\
\simeq&\delta^!\RintHom(q^{-1}_2( G_ 1\boxtimes \cdots \boxtimes G_\ell), q_1^! Ri_{\Delta \ast} F)\\
\simeq&\RintHom(\delta^{-1}q^{-1}_2(G_1\boxtimes \cdots \boxtimes G_\ell), \delta^! q_1^! Ri_{\Delta\ast} F)\\
\simeq&\RintHom(G_1\boxtimes \cdots \boxtimes G_\ell, Ri_{\Delta \ast} F)\\
\simeq& Ri_{\Delta \ast}\RintHom(G_1 \otimes \cdots \otimes G_\ell, F),
\end{align*}
where $\delta$ is a diagonal embedding $X^\ell \to (X^2)^\ell, (x_1, \dots, x_\ell)\mapsto (x_1, x_1, x_2,x_2,\dots, x_\ell,x_\ell)$. The first isomorphism follows from Proposition \ref{muiso}. The second isomorphism follows from Proposition \ref{homiso}.
\endproof

Next we extend the result in \cite[Corollary A.2]{U} to subanalytic sheaf theory.
\begin{prop}\label{duchida}
{\it
Let $\tau:E \to Z$ be a vector bundle and $E^\ast$ be its dual. Set $P = \left\{(x,y)\in E \times_Z E^\ast ; \langle x, y \rangle \geq 0\right\}$. For a conic object $G\in D^b_{\mathbb{R}^+}(k_{E_{sa}})$, there is a natural isomorphism
\[
\tau^{-1}R\tau_\ast G \otimes \omega_{E/Z} \simeq Rp_{1!!}p_2^{-1}G^\wedge
\]
and a commutative diagram
% https://q.uiver.app/?q=WzAsNCxbMCwwLCJHIl0sWzAsMSwiXFx0YXVeey0xfVxcdGF1X1xcYXN0IEciXSxbMSwxLCJScF97MSF9cF8yXnstMX1HXlxcd2VkZ2UiXSxbMSwwLCJScF97MSF9KHBfMl57LTF9R15cXHdlZGdlKV9QIl0sWzEsMF0sWzIsMV0sWzIsM10sWzMsMF1d
\[
\begin{tikzcd}
	G & Rp_{1!!}(p_2^{-1}G^\wedge)_P \otimes \omega_{E/Z} \arrow[l]\\
	\tau^{-1}R\tau_\ast G \arrow[u] & Rp_{1!!}p_2^{-1}G^\wedge \otimes \omega_{E/Z}, \arrow[l] \arrow[u] 
\end{tikzcd}
\]

where every horizontal arrow is an isomorphism.
}
\end{prop}

\proof
We have a commutative diagram
% https://q.uiver.app/?q=WzAsNixbMCwwLCJHIl0sWzAsMSwiXFx0YXVeey0xfVxcdGF1X1xcYXN0IEciXSxbMSwxLCJcXHRhdV57LTF9XFx0YXVfXFxhc3QgUnBfezEhfShwXzJeey0xfUdeXFx3ZWRnZSlfUCJdLFsxLDAsIlJwX3sxIX0ocF8yXnstMX1HXlxcd2VkZ2UpX1AiXSxbMiwwLCJScF97MSF9cF8yXnstMX1HXlxcd2VkZ2UiXSxbMiwxLCJcXHRhdV57LTF9XFx0YXVfXFxhc3QgUnBfezEhfXBfMl57LTF9R15cXHdlZGdlIl0sWzEsMF0sWzIsMSwiXFxzaW0iLDJdLFsyLDNdLFszLDAsIlxcc2ltIiwyXSxbNCwzXSxbNSwyLCJcXGFscGhhIiwyXSxbNSw0LCJcXGJldGEiLDJdXQ==
\[
\begin{tikzcd}
	G & {Rp_{1!!}(p_2^{!}G^\wedge)_P} \arrow[l, "\sim"'] & {Rp_{1!!}p_2^{!}G^\wedge} \arrow[l] \\
	{\tau^{-1}R\tau_\ast G} \arrow[u] & {\tau^{-1}R\tau_\ast Rp_{1!!}(p_2^{!}G^\wedge)_P} \arrow[l, "\sim"'] \arrow[u] & {\tau^{-1}R\tau_\ast Rp_{1!!}p_2^{!}G^\wedge,} \arrow[l, "\alpha"'] \arrow[u, "\beta"']
\end{tikzcd}
\]

where each horizontal arrow of the left square is an isomorphism by \cite[Theorem 3.2.6]{Lu2}.
To show $\alpha$ is an isomorphism, it is enough to prove
\[R\tau_\ast Rp_{1!!}H \to R\tau_\ast Rp_{1!!}H_P\]
is an isomorphism for a conic object $H \in D^b_{\mathbb{R}^+}(k_{{E\times_M E^\ast}_{sa}})$ on ${E\times_M E^\ast}_{sa}$. We make the use of the following Cartesian diagram:
% https://q.uiver.app/?q=WzAsNyxbMCwwLCJFXFx0aW1lc19NIEVeXFxhc3QiXSxbMSwwLCJNXFx0aW1lcyBFXlxcYXN0Il0sWzAsMSwiRSJdLFsxLDEsIk0iXSxbMiwwLCJFXFx0aW1lc19NIEVeXFxhc3QiXSxbMiwxLCJFIl0sWzMsMCwiRV5cXGFzdCJdLFswLDEsIlxcdGF1XlxccHJpbWUiXSxbMiwzLCJcXHRhdSJdLFsxLDMsInBfMV5cXHByaW1lIiwyXSxbMCwyLCJwXzEiLDJdLFsxLDQsImleXFxwcmltZSIsMCx7InN0eWxlIjp7InRhaWwiOnsibmFtZSI6Imhvb2siLCJzaWRlIjoidG9wIn19fV0sWzMsNSwiaSIsMCx7InN0eWxlIjp7InRhaWwiOnsibmFtZSI6Imhvb2siLCJzaWRlIjoidG9wIn19fV0sWzQsNSwicF8xIiwyXSxbNCw2LCJwXzIiXV0=
\[
\begin{tikzcd}
	{E\times_Z E^\ast} \arrow[r, "{\tau^\prime}"] \arrow[d, "{p_1}"'] & {Z\times E^\ast} \arrow[r, hook, "{i^\prime}"] \arrow[d, "{p_1^\prime}"'] & {E\times_Z E^\ast} \arrow[r, "{p_2}"] \arrow[d, "{p_1}"'] & {E^\ast} \\
	E \arrow[r, "\tau"'] & Z \arrow[r, hook, "i"'] & E. &
\end{tikzcd}
\]
%By Lemma \ref{lempd},  Proposition 2.5.2 and Proposition 2.5.3 of \cite{Lu2} we have $Rp_{1!!}H_P \in D^b_{\mathbb{R}^+}(k_{E_{sa}})$. 
We have the commutative diagram,
% https://q.uiver.app/?q=WzAsNixbMCwxLCJcXHRhdV9cXGFzdCBScF97MSF9SCJdLFswLDAsIiBcXHRhdV9cXGFzdCBScF97MSF9SF9QIl0sWzEsMCwiaV57LTF9UnBfezFcXGFzdH1IX1AiXSxbMSwxLCJpXnstMX1ScF97MVxcYXN0fUgiXSxbMiwwLCJScF5cXHByaW1lX3sxXFxhc3R9IGlee1xccHJpbWUtMX1IX1AiXSxbMiwxLCJScF5cXHByaW1lX3sxXFxhc3R9IGlee1xccHJpbWUtMX1IIl0sWzAsMV0sWzEsMiwiXFxzaW0iXSxbMCwzLCJcXHNpbSJdLFsyLDQsIlxcc2ltIl0sWzMsNSwiXFxzaW0iXSxbMywyXSxbNSw0XV0=
\[
\begin{tikzcd}
	{R\tau_\ast Rp_{1!!}H_P} \arrow[r,"\sim"] & {i^{-1}Rp_{1!!}H_P} \arrow[r, "\sim"] & {Rp^\prime_{1!!} i^{\prime-1}H_P} \\
	{R\tau_\ast Rp_{1!!}H} \arrow[r,"\sim"] \arrow[u] & {i^{-1}Rp_{1!!}H} \arrow[r, "\sim"] \arrow[u] & {Rp^\prime_{1!!} i^{\prime-1}H,} \arrow[u]
\end{tikzcd}
\]
where horizontal isomorphisms follows from
$R\tau_\ast  \simeq i^{-1}$ ( \cite[Lemma 3.1.5 ]{Lu2})
and the base change formula (\cite[Proposition 2.2.8]{Lu1}).
Since
\[i^\prime(Z\times_Z E^\ast) \subset P\]
third vertical arrow is an isomorphism. Thus $\alpha$ is an isomorphism.
We shall prove $\beta$ is an isomorphism. Since $p_2 = p_2 \circ i^\prime \circ \tau^\prime$ and by the base change formula(\cite[Proposition 2.2.9]{Lu1}), we have
\begin{align*}
 Rp_{1!!} p^{-1}_2 G^\wedge \otimes \omega_{E/Z} 
\simeq& Rp_{1!!} \tau^{\prime-1} i^{\prime-1}p^{-1}_2G^\wedge \otimes \omega_{E/Z}\\
\simeq& \tau^{-1} i^{-1}Rp_{1!!}p^{-1}_2G^\wedge \otimes \omega_{E/Z}\\
\simeq& \tau^{-1}R\tau_\ast Rp_{1!!}p_2^{!}G^\wedge.
\end{align*}
\endproof
\begin{prop}\label{fiso}
{\it
Let $G$ be a multi-conic object on $E^\ast_{sa}$. We have
\[G^{\lor^\ast_E}  \simeq Rp^{\prime\prime}_{1!!} p^{\prime\prime-1}_2 G \otimes \omega_{E/Z}\simeq Rp_{1!!} (p_2^! G)_P.\]
}
(For the notation see Subsection \ref{Sec: multimicro}.)
\end{prop}
\proof
By induction assume $\ell=2$. We have the following diagram:
% https://q.uiver.app/?q=WzAsNixbMSwwLCJQIl0sWzIsMCwiUF8xIFxcdW5kZXJzZXR7Wn17XFx0aW1lc30gRV8yXlxcYXN0Il0sWzMsMCwiRV5cXGFzdCJdLFsyLDEsIkVfMSBcXHVuZGVyc2V0e1p9e1xcdGltZXN9IEVeXFxhc3RfMiJdLFsxLDEsIkVfMSBcXHVuZGVyc2V0e1p9e1xcdGltZXN9IFBfMiJdLFswLDEsIkUiXSxbMCwxLCJwXntcXHByaW1lXFxwcmltZX1fezIsMn0iLDJdLFsxLDMsInBee1xccHJpbWVcXHByaW1lfV97MSwxfSJdLFswLDQsInBee1xccHJpbWVcXHByaW1lfV97MSwxfSJdLFs0LDMsInBee1xccHJpbWVcXHByaW1lfV97MiwyfSIsMl0sWzEsMiwicF57XFxwcmltZVxccHJpbWV9X3sxLDJ9IiwyXSxbNCw1LCJwXntcXHByaW1lXFxwcmltZX1fezIsMX0iXSxbMCw1LCJwXntcXHByaW1lXFxwcmltZX1fMSIsMix7ImN1cnZlIjoxfV0sWzAsMiwicF57XFxwcmltZVxccHJpbWV9XzIiLDAseyJjdXJ2ZSI6LTN9XSxbMCwxLCJ7fVxcIFxcIFxcIFxcIFxcc3F1YXJlIiwyLHsib2Zmc2V0Ijo1LCJzaG9ydGVuIjp7InNvdXJjZSI6MjB9LCJzdHlsZSI6eyJib2R5Ijp7Im5hbWUiOiJub25lIn0sImhlYWQiOnsibmFtZSI6Im5vbmUifX19XV0=

\[
\begin{tikzcd}
	& P \arrow[r, "{p^{\prime\prime}_{2,2}}"] \arrow[d, "{p^{\prime\prime}_{1,1}}"'] \arrow[ld, "{p^{\prime\prime}_1}"', bend right] \arrow[rr, bend left, "{p^{\prime\prime}_2}"] \ar[dr, phantom, "\square"]& {P_1 \underset{Z}{\times} E_2^\ast} \arrow[d, "{p^{\prime\prime}_{1,1}}"] \arrow[r, "{p^{\prime\prime}_{1,2}}"'] & {E^\ast} \\
	{E}  & {E_1 \underset{Z}{\times} P_2} \arrow[l, "{p^{\prime\prime}_{2,1}}"'] \arrow[r, "{p^{\prime\prime}_{2,2}}"'] & {E_1 \underset{Z}{\times} E^\ast_2.} &
\end{tikzcd}
\]
Then by Proposition \ref{embiso} we have 
\begin{align*}
G^{\lor^\ast_1 \lor^\ast_2} &\simeq  Rp^{}_{2,1!!} (p_{2,2}^{-1} Rp^{}_{1,1!!} (p^{-1}_{1,2} G)_{P_1})_{P_2} \otimes \omega_{E/Z}\\
&\simeq Rp^{\prime\prime}_{2,1!!} p_{2,2}^{\prime\prime-1} Rp^{\prime\prime}_{1,1!!} p_{1,2}^{\prime\prime-1}G \otimes  \omega_{E/Z}\\
&\simeq Rp^{\prime\prime}_{2,1!!} Rp^{\prime\prime}_{1,1!!} p_{2,2}^{\prime\prime-1} p_{1,2}^{\prime\prime-1}G \otimes \omega_{E/Z} \\
 &\simeq Rp^{\prime\prime}_{1!!} p^{\prime\prime -1}_2 G \otimes \omega_{E/Z}.
\end{align*}
Assume $\ell>2$. Set $E^\prime := \underset{Z, 2 \leq i \leq \ell}{\times}{E_i}$, $E^{\prime\ast}:= \underset{Z, 2 \leq i \leq \ell}{\times}{E^{\prime\ast}_i}$, $P_{E^\prime}:= \underset{Z, 2 \leq i \leq \ell}{\times}{P_i}$. We have the following commutative diagram:
% https://q.uiver.app/?q=WzAsNixbMSwwLCJQIl0sWzIsMCwiUF8xIFxcdW5kZXJzZXR7Wn17XFx0aW1lc30gRV57XFxwcmltZVxcYXN0fSJdLFszLDAsIkVeXFxhc3QiXSxbMiwxLCJFXzEgXFx1bmRlcnNldHtafXtcXHRpbWVzfSBFXntcXHByaW1lXFxhc3R9Il0sWzEsMSwiRV8xIFxcdW5kZXJzZXR7Wn17XFx0aW1lc30gUF97RV5cXHByaW1lfSJdLFswLDEsIkUiXSxbMCwxLCJwXntcXHByaW1lXFxwcmltZX1fe0VeXFxwcmltZSwyfSIsMl0sWzEsMywicF57XFxwcmltZVxccHJpbWV9X3sxLDF9Il0sWzAsNCwicF57XFxwcmltZVxccHJpbWV9X3sxLDF9Il0sWzQsMywicF57XFxwcmltZVxccHJpbWV9X3tFXlxccHJpbWUsMn0iLDJdLFsxLDIsInBee1xccHJpbWVcXHByaW1lfV97MSwyfSIsMl0sWzQsNSwicF57XFxwcmltZVxccHJpbWV9X3tFXlxccHJpbWUsMX0iXSxbMCw1LCJwXntcXHByaW1lXFxwcmltZX1fMSIsMix7ImN1cnZlIjoxfV0sWzAsMiwicF57XFxwcmltZVxccHJpbWV9XzIiLDAseyJjdXJ2ZSI6LTJ9XSxbMCwxLCJ7fVxcIFxcIFxcIFxcIFxcc3F1YXJlIiwyLHsib2Zmc2V0Ijo1LCJzaG9ydGVuIjp7InNvdXJjZSI6MjB9LCJzdHlsZSI6eyJib2R5Ijp7Im5hbWUiOiJub25lIn0sImhlYWQiOnsibmFtZSI6Im5vbmUifX19XV0=
\[
\begin{tikzcd}
	& P \arrow[r, "{p^{\prime\prime}_{E^\prime,2}}"] \arrow[d, "{p^{\prime\prime}_{1,1}}"'] \arrow[ld, "{p^{\prime\prime}_1}"', bend right] \arrow[rr, bend left, "{p^{\prime\prime}_2}"] \ar[dr, phantom, "\square"] & {P_1 \underset{Z}{\times} E^{\prime\ast}} \arrow[d, "{p^{\prime\prime}_{1,1}}"] \arrow[r, "{p^{\prime\prime}_{1,2}}"'] & {E^\ast} \\
	{E} & {E_1 \underset{Z}{\times} P_{E^\prime}} \arrow[r, "{p^{\prime\prime}_{E^\prime,2}}"'] \arrow[l, "{p^{\prime\prime}_{E^\prime,1}}"'] & {E_1 \underset{Z}{\times} E^{\prime\ast}.} &
\end{tikzcd}
\]
By induction hypothesis, we have
\begin{align*}
(G^{\lor^\ast_1})^{\lor^\ast_{E^\prime}} \simeq Rp^{\prime\prime}_{1!!} p^{\prime\prime -1}_2 G \otimes \omega_{E/Z}.
\end{align*}
Therefore, the induction proceeds.
\endproof
We need to extend the result \cite[Proposition 3.2.9]{Lu2}:
\begin{lem}\label{lempd}
{\it
Let $E_1, E_2$ be vector bundles over $Z$, a biconic object $F \in D^b_{\mathbb{R}^+ \times \mathbb{R}^+}(k_{{E_1 \times_Z E_2 }_{sa}})$ on  ${E_1 \times_Z E_2 }_{sa}$ and $p: E_1 \times_Z E_2 \to E_1$
then $Rp_{!!} F \in D^b_{\mathbb{R}^+}(k_{E_{1sa}})$. Here the first $\mathbb{R}^+$-action on ${E_1 \times_Z E_2 }$ is either $((z; \xi_1, \xi_2), t_1) \mapsto (z; t_1\xi_1, t_1\xi_2)$ or $((z; \xi_1, \xi_2), t_1) \mapsto (z; t_1\xi_1, \xi_2)$ and the second one
%$({E_1 \times_Z E_2 })\times\mathbb{R}^+\to E_1 \times_Z E_2$ on ${E_1 \times_Z E_2 }$
 is $( (z; \xi_1, \xi_2),t_2) \mapsto (z; \xi_1, t_2\xi_2)$.
Moreover,
\[
(Rp_{!!} F)^{\wedge_{E_1}} \simeq Rp^\prime_{!!} F^{\wedge_{E_1}},
\]
where $p^\prime: E_1^\ast \times_Z E_2 \to E_1^\ast$.
}
\end{lem}
\proof
The natural embedding $i: E_1 \to E_1 \times_Z E_2$ is a conic morphism with respect to the first action on ${E_1 \times_Z E_2 }$. By  \cite[Proposition 2.5.13]{Lu2} we get $i^{!} F\in D^b_{\mathbb{R}^+}(k_{E_{1sa}})$.  Since $F$ is a conic object on ${E_1 \times_Z E_2}_{sa}$ with respect to the second action, by \cite[Lemma 3.1.5]{Lu2} we obtain $ i^{!} F\simeq Rp_{!!}F$. Therefore $Rp_{!!}F\in D^b_{\mathbb{R}^+}(k_{E_{1sa}})$. Consider the diagram below where the squares are Cartesian:
% https://q.uiver.app/#q=WzAsNyxbMCwwLCJFXlxcYXN0XzEgXFx1bmRlcnNldHtafXtcXHRpbWVzfSBFXzIiXSxbMSwwLCJFXlxcYXN0XzEiXSxbMSwxLCJQXlxccHJpbWUiXSxbMSwyLCJFXzEiXSxbMCwxLCJQXlxccHJpbWUgXFx1bmRlcnNldHtafXtcXHRpbWVzfSBFXzIiXSxbMCwyLCJFXzEgXFx1bmRlcnNldHtafXtcXHRpbWVzfSBFXzIiXSxbMiwxLCJFXzFcXHVuZGVyc2V0e1p9e1xcdGltZXN9IEVeXFxhc3RfMSJdLFswLDEsInBeXFxwcmltZSJdLFsyLDEsInFfMiIsMl0sWzIsMywicV8xIl0sWzQsMiwiXFx0aWxkZXtwfSJdLFs0LDAsInFeXFxwcmltZV8yIl0sWzQsNSwicV8xXlxccHJpbWUiLDJdLFs1LDMsInAiXSxbNCwyLCJcXHNxdWFyZSIsMix7Im9mZnNldCI6NSwic3R5bGUiOnsiYm9keSI6eyJuYW1lIjoibm9uZSJ9LCJoZWFkIjp7Im5hbWUiOiJub25lIn19fV0sWzAsMSwiXFxzcXVhcmUiLDIseyJvZmZzZXQiOjUsInN0eWxlIjp7ImJvZHkiOnsibmFtZSI6Im5vbmUifSwiaGVhZCI6eyJuYW1lIjoibm9uZSJ9fX1dLFs2LDEsInBfMiIsMl0sWzYsMywicF8xIl0sWzIsNiwiaiIsMix7InN0eWxlIjp7InRhaWwiOnsibmFtZSI6Imhvb2siLCJzaWRlIjoidG9wIn19fV1d
\[
\begin{tikzcd}
	{E^\ast_1 \underset{Z}{\times} E_2} \arrow[r, "p^\prime"]\ar[dr, phantom, "\square"]  & {E^\ast_1} \\
	{P^\prime \underset{Z}{\times} E_2} \arrow[u, "q^\prime_2"] \arrow[d, "q_1^\prime"'] \arrow[r, "\tilde{p}"]\ar[dr, phantom, "\square"]  & {P^\prime} \arrow[u, "q_2"] \arrow[d, "q_1"] \arrow[r, hook, "j"'] & {E_1 \underset{Z}{\times} E^\ast_1} \arrow[lu, "p_2"] \arrow[ld, "p_1"] \\
	{E_1 \underset{Z}{\times} E_2} \arrow[r, "p"] & {E_1,}
\end{tikzcd}
\]
here $P^\prime  := \{(\eta, \xi)\in E_1 \underset{Z}{\times} E^\ast_1; \langle \eta, \xi \rangle \leq 0\}$.
Notice that $Rq_{2!!} q^{-1}_1 G \simeq G^{\wedge_{E_1}}$ for a conic object $G$ on $E_{1sa}$. In fact,
\begin{align*}
Rp_{2!!}( p^{-1}_1 G)_{P^\prime}&\simeq Rp_{2!!} Rj_{!!} j^{-1} p^{-1}_1 G\\
&\simeq Rq_{2!!} q^{-1}_1 G.
\end{align*}

Therefore:
\begin{align*}
(Rp_{!!} F)^{\wedge_{E_1}} &\simeq Rq_{2!!}q_{1}^{-1} Rp_{!!} F\\
&\simeq Rq_{2!!}R\tilde{p}_{!!} q^{\prime -1}_1F\\
&\simeq Rp_{!!}^\prime Rq^\prime_{2!!}q^{\prime-1}_1 F\\
&\simeq Rp_{!!}^\prime( F^{\wedge_{E_1}}).
\end{align*}
\endproof
\begin{prop}\label{iso}
{\it
Let $F$ be a multi-conic object on $E_{sa}$. We have a natural isomorphism
\[\tau^{-1}R\tau_\ast F  \simeq Rp_{1!!}p_2^{!}F^{\wedge_E}\]
and the natural morphism $\tau^{-1}R\tau_\ast F \to F$ is embedded to the following distinguished triangle:
% https://q.uiver.app/#q=WzAsNCxbMSwwLCJcXHRhdV57LTF9XFx0YXVfXFxhc3QgRiJdLFsyLDAsIkYiXSxbMywwLCJ7fSJdLFswLDAsIlJwX3sxISF9KHBfezJ9XnshfUZee1xcd2VkZ2VfRX0pX3tQXnstfX0iXSxbMCwxXSxbMSwyLCIrMSJdLFszLDBdXQ==
\[\begin{tikzcd}
	{Rp_{1!!}(p_{2}^{!}F^{\wedge_E})_{P^{-}}} & {\tau^{-1}R\tau_\ast F} & F & {{}.}
	\arrow[from=1-2, to=1-3]
	\arrow["{+1}", from=1-3, to=1-4]
	\arrow[from=1-1, to=1-2]
\end{tikzcd}\]
}
(For the notation see Subsection \ref{Sec: multimicro}.)
\end{prop}
\proof
If $\ell =1$, the result follows by the Proposition \ref{duchida}. Assume $\ell >1$, and $E^\prime := \varprod_{Z,i=2}^\ell E_i$, $E^{\prime\ast}:=\varprod_{Z, i=2}^\ell E^\ast_i$, and $P^\prime_{E^\prime}:= \varprod_{Z, i=2}^\ell P^\prime_i$. Moreover, let $\wedge_{E^\prime}$(resp. $\lor^\ast_{E^\prime}$) be the composition of $\wedge_i$(resp. $\lor^\ast_{i}$) for $i =2, \dots, \ell$. Consider the following diagram:
% https://q.uiver.app/?q=WzAsNixbMSwwLCJFXFx1bmRlcnNldHtafXtcXHRpbWVzfUVeXFxhc3QiXSxbMiwwLCJFXzEgXFx1bmRlcnNldHtafXtcXHRpbWVzfSBFXlxcYXN0Il0sWzIsMSwiRV8xIFxcdW5kZXJzZXR7Wn17XFx0aW1lc31FXntcXHByaW1lXFxhc3R9Il0sWzEsMSwiRV8xIFxcdW5kZXJzZXR7Wn17XFx0aW1lc30gRV5cXHByaW1lIFxcdW5kZXJzZXR7Wn17XFx0aW1lcyBFXntcXHByaW1lXFxhc3R9fSJdLFswLDEsIkUiXSxbMywwLCJFXlxcYXN0Il0sWzAsMSwiXFxzcXVhcmUiLDIseyJvZmZzZXQiOjQsInN0eWxlIjp7ImJvZHkiOnsibmFtZSI6Im5vbmUifSwiaGVhZCI6eyJuYW1lIjoibm9uZSJ9fX1dLFsxLDIsInBfezEsMX0iXSxbMCwzLCJwX3sxLDF9IiwyXSxbMywyLCJwX3tFXlxccHJpbWUsMn0iLDJdLFswLDQsInBfMSIsMix7ImN1cnZlIjoyfV0sWzMsNCwicF97RV5cXHByaW1lLDF9IiwyXSxbMSw1LCJwX3sxLDJ9IiwyXSxbMCw1LCJwXzIiLDAseyJjdXJ2ZSI6LTN9XSxbMCwxLCJwX3tFXlxccHJpbWUsMn0iLDJdXQ==
\[\begin{tikzcd}
    & {E\underset{Z}{\times}E^\ast} \arrow[r, "{p_{E^\prime,2}}"] \arrow[d, "{p_{1,1}}"'] \arrow[ld, "{p_1}"', bend right=12pt] \arrow[rr, "{p_2}", bend left=18pt] \ar[dr, phantom, "\square"]& {E_1 \underset{Z}{\times} E^\ast} \arrow[r, "{p_{1,2}}"'] \arrow[d, "{p_{1,1}}"] & {E^\ast} \\
    E \arrow[r, "{p_{E^\prime,1}}"'] & {E_1 \underset{Z}{\times} E^\prime \underset{Z}{\times} E^{\prime\ast}} \arrow[r, "{p_{E^\prime,2}}"'] & {E_1 \underset{Z}{\times} E^{\prime\ast}.} &
\end{tikzcd}\]
By the commutative diagram
% https://q.uiver.app/?q=WzAsNCxbMCwwLCJFIl0sWzEsMCwiRV5cXHByaW1lIl0sWzAsMSwiRV8xIl0sWzEsMSwiWiJdLFswLDEsIlxcdGF1XzEiXSxbMCwyLCJcXHRhdV97RV5cXHByaW1lfSJdLFsxLDMsIlxcdGF1X3tFXlxccHJpbWV9Il0sWzIsMywiXFx0YXVfMSJdLFswLDEsIlxcc3F1YXJlIiwyLHsib2Zmc2V0Ijo1LCJzdHlsZSI6eyJib2R5Ijp7Im5hbWUiOiJub25lIn0sImhlYWQiOnsibmFtZSI6Im5vbmUifX19XV0=
\[
\begin{tikzcd}
	E \arrow[r, "\tau_1"] \arrow[d, "\tau_{E^\prime}"'] & {E^\prime} \arrow[d, "\tau_{E^\prime}"] \\
	{E_1} \arrow[r, "\tau_1"'] & Z
	\arrow[draw=none, from=1-1, to=2-2, "\square", phantom]
\end{tikzcd}\]
and \cite[Proposition 2.4.8]{Lu1} and the dual base change formula we have
\begin{align*}
\tau^{-1}_1 R\tau_{1\ast}\tau^{-1}_{E^\prime}R\tau_{E^\prime\ast}  F &\simeq \tau^{-1}_1 \tau^{-1}_{E^\prime} R\tau_{1\ast}  R\tau_{E^\prime\ast}  F\\
&\simeq \tau^{-1}R\tau_\ast F.
\end{align*}
Hence, by Proposition \ref{duchida} and induction hypothesis, we obtain the following commutative diagram:
% https://q.uiver.app/#q=WzAsOSxbMCwwLCJGIl0sWzEsMCwiRl57XFx3ZWRnZV97RV5cXHByaW1lfXtcXGxvcl97RV5cXHByaW1lfX1eXFxhc3R9Il0sWzIsMCwiRl57XFx3ZWRnZV97RV5cXHByaW1lfVxcbG9yX3tFXlxccHJpbWV9XlxcYXN0XFx3ZWRnZV8xe1xcbG9yXzF9XlxcYXN0fSJdLFswLDEsIlxcdGF1XnstMX1fe0VeXFxwcmltZX1SXFx0YXVfe0VeXFxwcmltZVxcYXN0fUYiXSxbMCwyLCJcXHRhdV57LTF9XzFSXFx0YXVfezFcXGFzdH0gXFx0YXVeey0xfV97RV5cXHByaW1lfVJcXHRhdV97RV5cXHByaW1lXFxhc3R9RiJdLFswLDMsIlxcdGF1XnstMX1SXFx0YXVfXFxhc3QgRiJdLFsyLDEsIlJwX3tFXlxccHJpbWUsMSEhfXBeeyF9X3tFXlxccHJpbWUsMn0oRl57XFx3ZWRnZV97RV5cXHByaW1lfX0pIl0sWzIsMiwiXFx0YXVeey0xfV8xUlxcdGF1X3sxXFxhc3R9UnBfe0VeXFxwcmltZSwxISF9cF57IX1fe0VeXFxwcmltZSwyfShGXntcXHdlZGdlX3tFXlxccHJpbWV9fSkiXSxbMiwzLCJScF97MSwxISF9cF57IX1fezEsMn0oKFJwX3tFXlxccHJpbWUsMSEhfXBeeyF9X3tFXlxccHJpbWUsMn0oRl57XFx3ZWRnZV97RV5cXHByaW1lfX0pKV57XFx3ZWRnZV8xfSkuIl0sWzEsMCwiXFxzaW0iLDJdLFsyLDEsIlxcc2ltIiwyXSxbMywwXSxbNCwzXSxbNSw0LCIiLDIseyJvZmZzZXQiOjEsInN0eWxlIjp7ImhlYWQiOnsibmFtZSI6Im5vbmUifX19XSxbNiwzLCJcXHNpbSIsMl0sWzcsNCwiXFxzaW0iLDJdLFs4LDUsIlxcc2ltIiwyXSxbNyw2XSxbNiwyXSxbOCw3LCJcXHNpbSJdLFs1LDQsIiIsMix7InN0eWxlIjp7ImhlYWQiOnsibmFtZSI6Im5vbmUifX19XV0=
\[\begin{tikzcd}
	F & {F^{\wedge_{E^\prime}{\lor_{E^\prime}}^\ast}} \arrow[l, "\sim"'] & {F^{\wedge_{E^\prime}\lor_{E^\prime}^\ast\wedge_1{\lor_1}^\ast}} \arrow[l, "\sim"'] \\
	{\tau^{-1}_{E^\prime}R\tau_{E^\prime\ast}F} \arrow[u] && {Rp_{E^\prime,1!!}p^{!}_{E^\prime,2}(F^{\wedge_{E^\prime}})} \arrow[ll, "\sim"'] \arrow[u] \\
	{\tau^{-1}_1R\tau_{1\ast} \tau^{-1}_{E^\prime}R\tau_{E^\prime\ast}F} \arrow[u] && {\tau^{-1}_1R\tau_{1\ast}Rp_{E^\prime,1!!}p^{!}_{E^\prime,2}(F^{\wedge_{E^\prime}})} \arrow[ll, "\sim"'] \arrow[u] \\
	{\tau^{-1}R\tau_\ast F} \arrow[u, shift right, no head] \arrow[u, shift left, no head] && {Rp_{1,1!!}p^{!}_{1,2}((Rp_{E^\prime,1!!}p^{!}_{E^\prime,2}(F^{\wedge_{E^\prime}}))^{\wedge_1}).} \arrow[ll, "\sim"'] \arrow[u, "\sim"]
\end{tikzcd}\]
By Lemma \ref{lorlor} below for a multi-conic object $H$ on $E_1 \times_Z E^{\prime\ast}_{sa}$ we have
$H^{\lor_{E^\prime}^\ast{\lor_1}} \simeq H^{{\lor_1}\lor_{E^\prime}^\ast}$.
Therefore, we have
\begin{align*}
H^{\lor_{E^\prime}^\ast\wedge_1} &\simeq (H^{\lor_{E^\prime}^\ast\lor_1})^{id \times a} \otimes \omega^{\otimes-1}_{E_1/Z}\\
&\simeq (H^{\lor_1\lor_{E^\prime}^\ast})^{id \times a} \otimes \omega^{\otimes-1}_{E_1/Z}\\
 &\simeq H^{\wedge_1\lor_{E^\prime}^\ast}.
\end{align*}
Thus, we get the following chain of isomorphisms by Proposition \ref{fiso} and \cite[Proposition 2.2]{HPY}.
\begin{align*}
F^{\wedge_{E^\prime}\lor_{E^\prime}^\ast\wedge_1{\lor_1}^\ast} &\simeq F^{\wedge_{E^\prime}\wedge_1\lor_{E^\prime}^\ast{\lor_1}^\ast}\\
&\simeq Rp^{\prime\prime}_{1!!}p^{\prime\prime-1}_2 F^{\wedge_E} \otimes \omega_{E/Z}.
\end{align*} 
Lastly by Lemma \ref{lempd} and \cite[Proposition 3.2.9]{Lu2} we have
\begin{align}
&{}Rp_{1,1!!}p^{-1}_{1,2}\left(\left(Rp_{E^\prime,1!!}p_{E^\prime,2}^{-1}(F^{\wedge_{E^\prime}})\right)^{\wedge_1}\right)\nonumber\\
&\simeq Rp_{1,1!!}p^{-1}_{1,2}R p^\prime_{E^\prime,1 !!} p^{\prime-1}_{E^\prime,2}(F^{\wedge_E})\label{morphisms}\\
&\simeq Rp_{1,1!!}R p_{E^\prime,1 !!} p^{\prime-1}_{1,2}p^{\prime-1}_{E^\prime,2}(F^{\wedge_E})\label{bach}\\
&\simeq Rp_{1,1!!}R p_{E^\prime,1 !!} p_2^{-1}(F^{\wedge_E})\nonumber\\
&\simeq Rp_{1!!}p_2^{-1}(F^{\wedge_E})\nonumber,
\end{align}
where morphisms on (\ref{morphisms}) are $p^\prime_{E^\prime, 2}: E^\prime \times_Z E^\ast \to E^\ast$ and $p^\prime_{E^\prime,1}: E^\ast_1 \times_Z E^\prime \times_Z E^{\prime\ast} \to E^\ast_1 \times_Z E^\prime$. On (\ref{bach}) we used the base change formula(\cite[Proposition 2.2.9]{Lu1}) on the Cartesian square as shown below.
% https://q.uiver.app/#q=WzAsNCxbMSwwLCJFXlxccHJpbWUgXFx1bmRlcnNldHtafXtcXHRpbWVzfSBFXlxcYXN0Il0sWzEsMSwiRV5cXGFzdF8xIFxcdW5kZXJzZXR7Wn17XFx0aW1lc30gRV5cXHByaW1lIl0sWzAsMSwiRV8xIFxcdW5kZXJzZXR7Wn17XFx0aW1lc31FXlxcYXN0XzFcXHVuZGVyc2V0e1p9e1xcdGltZXN9RV5cXHByaW1lIl0sWzAsMCwiRSBcXHVuZGVyc2V0e1p9e1xcdGltZXN9IEVeXFxhc3QiXSxbMCwxLCJwXlxccHJpbWVfe0VeXFxwcmltZSwxfSJdLFsyLDEsInBfezEsMn0iLDJdLFszLDAsInBeXFxwcmltZV97MSwyfSJdLFszLDIsInBfe0VeXFxwcmltZSwxfSIsMl1d
\[\begin{tikzcd}
	{E \underset{Z}{\times} E^\ast} \arrow[r, "{p^\prime_{1,2}}"] \arrow[d, "{p_{E^\prime,1}}"'] & {E^\prime \underset{Z}{\times} E^\ast} \arrow[d, "{p^\prime_{E^\prime,1}}"] \\
	{E_1 \underset{Z}{\times} E^\ast_1 \underset{Z}{\times} E^\prime} \arrow[r, "{p_{1,2}}"'] & {E^\ast_1 \underset{Z}{\times} E^\prime}
\end{tikzcd}\]
By the above isomorphisms,
% https://q.uiver.app/#q=WzAsMTIsWzEsMCwiXFx0YXVeey0xfVJcXHRhdV9cXGFzdCBGIl0sWzIsMCwiRiJdLFsxLDEsIlJwX3sxLDEhIX1wXnshfV97MSwyfShScF97RV5cXHByaW1lLDEhIX1wXiFfe0VeXFxwcmltZSwyfUZee1xcd2VkZ2Vfe0VeXFxwcmltZX19KV57XFx3ZWRnZV8xfSJdLFsyLDEsIkZee1xcd2VkZ2Vfe0VeXFxwcmltZX1cXGxvcl97RV5cXHByaW1lfV5cXGFzdFxcd2VkZ2VfMXtcXGxvcl8xfV5cXGFzdH0iXSxbMiwyLCJScF57XFxwcmltZVxccHJpbWV9X3sxISF9cF57XFxwcmltZVxccHJpbWUtMX1fMiAoRl57XFx3ZWRnZV9FfSlcXG90aW1lcyBcXG9tZWdhX3tFL1p9Il0sWzEsMiwiUnBfezEhIX1wXiFfMkZee1xcd2VkZ2VfRX0iXSxbMCwyLCJScF97MSEhfShwXnshfV8yKEZee1xcd2VkZ2Vfe0V9fSkpX3tQXnstfX0iXSxbMywyLCJ7fSJdLFswLDMsIlJwX3sxISF9KHBeeyF9XzIoRl57XFx3ZWRnZV97RX19KSlfe1Beey19fSJdLFsxLDMsIlJwX3sxISF9cF4hXzJGXntcXHdlZGdlX0V9Il0sWzIsMywiUnBfezEhIX0ocF4hXzJGXntcXHdlZGdlX0V9KV97UH0iXSxbMywzLCIuIl0sWzAsMV0sWzIsMCwiXFxzaW0iLDJdLFszLDEsIlxcc2ltIiwyXSxbMiwzXSxbNSw0XSxbNSwyLCJcXHNpbSIsMl0sWzQsMywiXFxzaW0iLDJdLFs2LDVdLFs0LDcsIisxIl0sWzYsOCwiIiwyLHsib2Zmc2V0IjotMSwic3R5bGUiOnsiaGVhZCI6eyJuYW1lIjoibm9uZSJ9fX1dLFs2LDgsIiIsMix7Im9mZnNldCI6MSwic3R5bGUiOnsiaGVhZCI6eyJuYW1lIjoibm9uZSJ9fX1dLFs4LDldLFs5LDEwXSxbNSw5LCIiLDIseyJzdHlsZSI6eyJoZWFkIjp7Im5hbWUiOiJub25lIn19fV0sWzUsOSwiIiwyLHsib2Zmc2V0IjoyLCJzdHlsZSI6eyJoZWFkIjp7Im5hbWUiOiJub25lIn19fV0sWzEwLDQsIlxcc2ltIiwyLHsic3R5bGUiOnsidGFpbCI6eyJuYW1lIjoiYXJyb3doZWFkIn0sImhlYWQiOnsibmFtZSI6Im5vbmUifX19XSxbMTAsMTEsIisxIl1d

\adjustbox{scale=0.8,center}
{
\begin{tikzcd}
	& {\tau^{-1}R\tau_\ast F} \arrow[r] & F \\
	& {Rp_{1,1!!}p^{!}_{1,2}(Rp_{E^\prime,1!!}p^!_{E^\prime,2}F^{\wedge_{E^\prime}})^{\wedge_1}} \arrow[r] \arrow[u, "\sim"'] & {F^{\wedge_{E^\prime}\lor_{E^\prime}^\ast\wedge_1{\lor_1}^\ast}} \arrow[u, "\sim"'] \\
	{Rp_{1!!}(p^{!}_2(F^{\wedge_{E}}))_{P^{-}}} \arrow[r]\arrow[d, no head, shift left] \arrow[d, no head, shift right=2] & {Rp_{1!!}p^!_2F^{\wedge_E}} \arrow[r] \arrow[u, "\sim"'] \arrow[d, no head, shift left] \arrow[d, no head, shift right=2] & {Rp^{\prime\prime}_{1!!}p^{\prime\prime-1}_2 (F^{\wedge_E})\otimes \omega_{E/Z}} \arrow[r, "{+1}"] \arrow[u, "\sim"'] \arrow[d, "\sim"] & {} \\
	{Rp_{1!!}(p^{!}_2(F^{\wedge_{E}}))_{P^{-}}} \arrow[r] & {Rp_{1!!}p^!_2F^{\wedge_E}} \arrow[r] & {Rp_{1!!}(p^!_2F^{\wedge_E})_{P}} \arrow[r, "{+1}"] & {.}
\end{tikzcd}
}

\endproof
\begin{lem}\label{lorlor}
\textit{
Let $E_1$ and $E^{\prime\ast}$ be as in the proof of Proposition \ref{iso}. For a multi-conic object $H$ on $E_1 \times_Z E^{\prime\ast}_{sa}$ we have
$H^{\lor_{E^\prime}^\ast{\lor_1}} \simeq H^{{\lor_1}\lor_{E^\prime}^\ast}$.
}
\end{lem}
\proof
Consider the following commutative diagram:
% https://q.uiver.app/#q=WzAsNixbMSwwLCJQXlxccHJpbWVfMSBcXHRpbWVzIFBfe0VeXFxwcmltZX1eXFxwcmltZSJdLFsyLDAsIlBeXFxwcmltZV8xIFxcdGltZXMgRV57XFxwcmltZVxcYXN0fSJdLFsyLDEsIkVeXFxhc3RfMVxcdGltZXMgRV57XFxwcmltZVxcYXN0fSJdLFsxLDEsIkVeXFxhc3RfMSBcXHRpbWVzX1ogUF5cXHByaW1lX3tFXlxccHJpbWV9Il0sWzAsMSwiRV5cXGFzdF8xIFxcdGltZXMgRV5cXHByaW1lIl0sWzMsMCwiRV8xXFx0aW1lcyBFXntcXHByaW1lXFxhc3R9Il0sWzAsMSwicF97RV5cXHByaW1lLDJ9XlxccHJpbWUiLDJdLFsxLDIsInBfezEsMn1eXFxwcmltZSJdLFswLDMsInBfezEsMn1eXFxwcmltZSJdLFszLDIsInBfe0VeXFxwcmltZSwgMn1eXFxwcmltZSIsMl0sWzMsNCwicF97RV5cXHByaW1lLCAxfV5cXHByaW1lIl0sWzAsNCwicF8xXlxccHJpbWUiLDIseyJjdXJ2ZSI6MX1dLFsxLDUsInBfezEsMX1eXFxwcmltZSJdLFswLDUsInBfMl5cXHByaW1lIiwwLHsiY3VydmUiOi0zfV1d
\[\begin{tikzcd}
	& {P^\prime_1 \underset{Z}{\times} P_{E^\prime}^\prime} \arrow[d, "{p_{1,2}^\prime}"] \arrow[dl, "{p_1^\prime}"'] \arrow[r, "{p_{E^\prime,2}^\prime}"] \arrow[rr, bend left=30, "{p_2^\prime}"] & {P^\prime_1 \underset{Z}{\times} E^{\prime\ast}} \arrow[d, "{p_{1,2}^\prime}"] \arrow[r, "{p_{1,1}^\prime}"] & {E_1 \underset{Z}{\times} E^{\prime\ast}} \\
	{E^\ast_1 \underset{Z}{\times} E^\prime} & {E^\ast_1 \underset{Z}{\times} P^\prime_{E^\prime}} \arrow[l, "{p_{E^\prime,1}^\prime}"] \arrow[r, "{p_{E^\prime,2}^\prime}"'] & {E^\ast_1 \underset{Z}{\times} E^{\prime\ast},}
\end{tikzcd}\]
where the square is Cartesian.
Then we have
\begin{align*}
H^{\lor_{E^\prime}^\ast{\lor_1}} \simeq Rp_{E^\prime,1\ast}^\prime p_{E^\prime,2}^{\prime!}Rp_{1,2\ast}^\prime p_{1,1}^{\prime!} H \simeq& Rp_{E^\prime,1\ast}^\prime Rp_{1 ,2\ast}^\prime p_{E^\prime,2}^{\prime!}p_{1,1}^{\prime!} H\\
\simeq& Rp_{1\ast}^\prime p^{\prime!}_{2}H.
\end{align*}
For the same reason, we obtain $H^{{\lor_1}\lor_{E^\prime}^\ast}\simeq Rp_{1\ast}^\prime p^{\prime!}_{2}H$.
\endproof

\begin{prop}\label{isom2}
\textit{
Let $F \in D^b(k_{X_{sa}})$. Let $i$ and $M$ be as in Proposition \ref{muiso}. Then:
\[R\pi_{!!} \mu^{sa}_\chi(F) \simeq R\Gamma_M(\mu^{sa}_\chi(F)) \simeq i^{-1}F \otimes \omega_{M/X}.\]
}
\end{prop}
\proof
Set $H:=\nu^{sa}_\chi(F)$ and $E:= S_\chi $, then by Proposition \ref{iso}
we have
\[
\tau^{-1}R\tau_\ast H \otimes \omega_{E/M} \simeq Rp_{1!!}p_2^{-1}H^{\wedge_E}.
\]
By the base change formula(\cite[Proposition 1.4.6]{Lu1}) on the Cartesian square,
\[\begin{tikzcd}
	{E \times_M E^\ast} \arrow[r, "{p_2}"] \arrow[d, "{p_1}"'] & {E^\ast} \arrow[d, "\pi"] \\
	E \arrow[r, "\tau"'] & M
\end{tikzcd}
\]
we have
\[\tau^{-1}R\pi_{!!} H^{\wedge_E} \simeq Rp_{1!!}p_2^{-1}H^{\wedge_E}.\]
Thus
\[\tau^{-1}R\pi_{!!} H^{\wedge_E} \simeq \tau^{-1}R\tau_\ast H \otimes \omega_{E/M}.\]
We apply $i^{-1}$ to get
\[R\pi_{!!} H^{\wedge_E} \simeq R\tau_\ast H \otimes \omega_{M/X}.\]
Since $R\tau_\ast H  \simeq F|_M$ (\cite[Proposition 2.6]{HPY}) we have,
\[ R\pi_{!!} \nu^{sa}_\chi(F)^{\wedge_E} \simeq F|_M \otimes \omega_{M/X}.\]
\endproof
\begin{prop}
\textit{
Let $G_1,\dots ,G_\ell \in D^b_{\mathbb{R}\text{-}c}(k_{X})$ and $F \in D^b(k_{X_{sa}})$.  We obtain the following isomorphisms:
\begin{align*}
R\pi_{!!}\mu hom^{sa}_{\widehat{\chi}}(G_1, G_2,\dots G_\ell,F) \simeq& \RintHom(G_1\boxtimes\dots  \boxtimes G_\ell,k_{X^\ell}) \otimes  Ri_{\Delta\ast}F\\
\simeq& Ri_{\Delta\ast}(D^\prime G_1\otimes \dots  \otimes D^\prime G_\ell \otimes  F),
\end{align*}
where we set $D^\prime F:= \RintHom(F, k_X)$ for $F \in D^b(k_{X_{sa}})$.
}
\end{prop}
\proof
We have the isomorphisms by Proposition \ref{isom2}:
\begin{align}
{}&R\pi_{!!}\mu hom^{sa}_{\widehat{\chi}}(G_1, G_2,\dots G_\ell,F) \nonumber\\
\simeq& \delta^{-1}\left(\RintHom(p^{-1}_2(G_1\boxtimes\dots  \boxtimes G_\ell), p^{!}_1 Ri_{\Delta\ast}F)\right)\otimes \omega_{X^\ell/\widehat{X}}\nonumber\\
%&\simeq&\delta^{-1}\left(\RintHom(p^{-1}_2(G_1\boxtimes\dots  \boxtimes G_\ell), p^{-1}_1 Ri_{\Delta\ast}F)\right)\otimes \delta^{-1}\omega_{(X^2)^\ell/X^\ell}\otimes \omega_{X^\ell/\widehat{X}}\nonumber\\
%&\simeq&\delta^{-1}\left(\RintHom(p^{-1}_2(G_1\boxtimes\dots  \boxtimes G_\ell), p^{-1}_1 Ri_{\Delta\ast}F)\right)\otimes \omega_{X^\ell/X^{\ell}}\nonumber\\%this is because KS90 (3.3.4)
\simeq& \delta^{-1}(\RintHom(G_1\boxtimes\dots  \boxtimes G_\ell,k_{X^\ell}) \boxtimes Ri_{\Delta\ast} F )\label{box}\\
\simeq& \RintHom(G_1\boxtimes\dots  \boxtimes G_\ell,k_{X^\ell}) \otimes  Ri_{\Delta\ast}F\nonumber\\
\simeq& Ri_{\Delta\ast}(D^\prime G_1\otimes \dots  \otimes D^\prime G_\ell \otimes  F)\label{4.7}.
\end{align}
Here $\delta: X^\ell \hookrightarrow (X^2)^\ell$ is a diagonal embedding. The isomorphism (\ref{box}) follows from \cite[Lemma 5.3.9]{Lu2}. (\ref{4.7}) follows from  \cite[Proposition 3.4.4]{KS90}.
\endproof
\begin{cor}[multi-microlocal Sato's triangle]
{\it
We get the distinguished triangle:
% https://q.uiver.app/?q=WzAsNCxbMCwwLCJGfF9NIFxcb3RpbWVzIFxcb21lZ2Ffe00vWH0iXSxbMSwwLCJSXFxHYW1tYV9NKEYpfF9NIl0sWzIsMCwiUlxcZG90XFxwaV9cXGFzdCBcXG11X1xcY2hpKEYpIl0sWzMsMCwie30iXSxbMCwxXSxbMSwyXSxbMiwzLCIrMSJdXQ==
\[\begin{tikzcd}
	{F|_M \otimes \omega_{M/X}} \arrow[r]& {R\Gamma_M(F)|_M} \arrow[r] & {R\dot\pi_\ast \mu^{sa}_\chi(F)} \arrow[r, "+1"] & {,} 
\end{tikzcd}\]
where $\dot\pi$ is the restriction of $\pi: S^\ast_\chi \to M$ to $S^\ast_\chi \setminus M$.
}
\end{cor}
\proof
We finish the proof by applying the above results(Proposition \ref{muiso}, Proposition \ref{isom2}) to the distinguished triangle:
% https://q.uiver.app/?q=WzAsNCxbMCwwLCJSXFxHYW1tYV9NKFxcbXVfXFxjaGkoRikpIl0sWzEsMCwiUlxccGlfXFxhc3QgXFxtdV9cXGNoaSAoRikiXSxbMiwwLCJSXFxkb3R7XFxwaX1fXFxhc3RcXG11X1xcY2hpKEYpIl0sWzMsMCwie30iXSxbMCwxXSxbMSwyXSxbMiwzLCIrMSJdXQ==
\[
\begin{tikzcd}
	{R\pi_{!!}\mu^{sa}_\chi(F)} \arrow[r] & {R\pi_\ast \mu^{sa}_\chi (F)} \arrow[r] & {R\dot{\pi}_\ast \mu^{sa}_\chi(F)} \arrow[r, "+1"] & {.} 
\end{tikzcd}
\]
\endproof
\begin{cor}
{\it
For $G_1, G_2,\dots,G_\ell \in D^b_{\mathbb{R}\text{-}c}(k_{X})$ and $F \in D^b(k_{X_{sa}})$, we have the following distinguished triangle:
% https://q.uiver.app/?q=WzAsNCxbMCwwLCJcXG92ZXJzZXR7bH17XFx1bmRlcnNldHtpPTF9e1xcb3RpbWVzfSB9XFxtYXRoYmZ7RH1HX2kgXFxvdGltZXMgIEYiXSxbMSwwLCJcXFJpbnRIb20oXFxvdmVyc2V0e2x9e1xcdW5kZXJzZXR7aT0xfXtcXG90aW1lc30gfUdfaSwgIEYpIl0sWzIsMCwiUlxcZG90e1xccGl9XypcXG11IGhvbV57c2F9X3tcXGhhdHtcXGNoaX19KEdfMSxcXGRvdHMsIEdfbDtGKSJdLFszLDAsInt9Il0sWzAsMV0sWzEsMl0sWzIsMywiKzEiXV0=
\[\begin{tikzcd}
	{\overset{\ell}{\underset{i=1}{\otimes} }D^\prime G_i \otimes F} \arrow[r] & {\RintHom(\overset{\ell}{\underset{i=1}{\otimes} }G_i, F)} \arrow[r] & {R\dot{\pi}_*\mu \text{hom}^{sa}_{\widehat{\chi}}(G_1, \dots, G_\ell; F)} \arrow[r, "+1"] & {,}
\end{tikzcd}\]
where $\dot{\pi}$ is the restriction of $\pi: T^\ast X \times_X \cdots \times_X T^\ast X \to X$ to $T^\ast X \times_X \cdots \times_X T^\ast X \setminus X$.
}
\end{cor}
\begin{ex}
Let $X=\mathbb{C}^2$ and $\widehat{\chi}= \{\Delta \times \mathbb{C}^2 , \mathbb{C}^2 \times \Delta \}$. Consider closed subsets
$M_1 = \{ y_1 = 0, x_1 \ge 0\} = \mathbb{R}_{\ge 0} \times  \mathbb{C}$ and
$M_2 = \{ y_2 = 0, x_2 \ge 0\} =  \mathbb{C}   \times \mathbb{R}_{\ge 0}$.

Since $D^\prime(\mathbb{C}_{M_1}) \simeq \mathbb{C}_{\{ y_1 = 0, x_1 > 0\}}[-1]$, $D^\prime(\mathbb{C}_{M_2}) \simeq \mathbb{C}_{\{ y_2 = 0, x_2 > 0\}}[-1]$ 
%and $\left(\mathbb{C}_{\{y_1 = 0, y_2 = 0, x_1>0, x_2>0\}}\right)_{(0,0)} =0$ 
we obtain
\[
\left(\RintHom(\mathbb{C}_{M_1\cap M_2}, R\rho_\ast \mathscr{O}_X)\right)_{(0,0)} \simeq \left(R\dot{\pi}_{\ast}\mu hom_{\widehat{\chi}}(\mathbb{C}_{M_1},\mathbb{C}_{M_2}; R\rho_\ast\mathscr{O}_X )\right)_{(0,0)}.
\]
\end{ex}

\section*{Acknowledgement}
The author is supported in part by JST SPRING, Grant Number JPMJSP2119.

\addcontentsline{toc}{section}{

\end{document}